\documentclass{article}
\usepackage{graphicx} 
\usepackage{amsfonts}
\usepackage[mathscr]{eucal}
\usepackage{mathrsfs}
\usepackage{amsmath}
\usepackage{amssymb}
\usepackage{amsmath}
\usepackage{amsthm}
\usepackage{mathrsfs}
\usepackage{todonotes}
\usepackage{enumitem}
\newtheorem{theorem}{Theorem}[section]
\newtheorem{corollary}{Corollary}
\newtheorem{lemma}[theorem]{Lemma}

\newtheorem{remark}{Remark}
\newtheorem{assumption}{Assumption}
\newtheorem{example}{Example}

\newcommand{\R}{\mathbb{R}}
\newcommand{\eps}{\varepsilon}
\newcommand{\Q}{\mathcal{Q}}
\DeclareMathOperator{\tr}{tr}
\DeclareMathOperator{\sym}{sym}
\def\dist{\mathrm{dist\,}}
\def\calA{\mathcal{A}}
\title{Derivation of the Reissner-Mindlin model from nonlinear elasticity}
\author{Tamara Fastovska\thanks{Humboldt-Universit\"at zu Berlin, Institut f\"ur Mathematik, Unter den Linden 6, 10099 Berlin, email: tamara.fastovska@hu-berlin.de\\
V.N. Karazin Kharkiv National University, Department of Mathematics and Computer Science, Svobody sq. 4, 61002 Kharkiv, email: fastovskaya@karazin.ua}, Janusz Ginster\thanks{Weierstrass Institute, Mohrenstrasse 39, 10117 Berlin, Germany, email: ginster@wias-berlin.de.} and Barbara Zwicknagl\thanks{Humboldt-Universit\"at zu Berlin, Institut f\"ur Mathematik, Unter den Linden 6, 10099 Berlin, email: barbara.zwicknagl@hu-berlin.de} }
\date{}

\begin{document}
\maketitle

\begin{abstract}
\noindent
We discuss how the Reissner-Mindlin plate model can be derived from three-dimensional finite elasticity in terms of $\Gamma$-convergence. The presence of transverse shear effects in the Reissner-Mindlin model requires to scale different components of the three-dimensional elastic strain differently. 
A main technical tool is then the combination of rigidity estimates for the deformation and suitably averaged versions.
\end{abstract}
\textbf{Keywords: }Reissner-Mindlin plate model, dimension reduction, finite elasticity\\
\textbf{Mathematics Subject Classification (2020)}: 49J45, 74B20, 74K20
\section{Introduction}
The rigorous justification of simplified 
models for the elastic behaviour of thin bodies has a long standing history. In this paper, we follow the approach to derive static lower-dimensional models from 
three-dimensional elasticity by means of $\Gamma$-convergence. The latter is a concept of convergence that -if complemented by appropriate compactness properties- ensures convergence of (almost) minimizers (see e.g. \cite{braides,DalMaso}). The idea to use $\Gamma$-convergence in the context of dimension reduction in elasticity dates back to \cite{ABP, FJM, LDR1,LDR2,Pan}, and since then, a large body of literature has been devoted to a refined analysis in various settings (for recent overviews see e.g. \cite{lewicka,mueller,PPG}), including for instance nonlinear bending models for shells as well as Kirchhoff-Love's or von K\'arm\'an models for plates
(see e.g. \cite{ADS,BGNPP,DBS,FJM,FGM_Der,GOW,LM,NV,BS} and the references therein) and models for beams or rods (see e.g. \cite{Davoli,EK,GG,MM1,MM2,Neu,Scardia,Scardia2} and the references therein).

\par 

In the present work we address the  question, whether it is possible to derive the Reissner-Mindlin plate model (\cite{M,R}) from nonlinear elasticity. Before discussing the related literature let us briefly explain the setting and the challenges. We denote by $S\subseteq \mathbb R^2$ the midsurface of a thin plate $S\times(-h/2,h/2)$, and consider the energy functionals for linearized and nonlinear Reissner-Mindlin models, namely
\begin{eqnarray}
\label{ener_mtl}
 J(\varphi, v)=\sum\limits_{i=1}^2a_i\int\limits_S (\varphi_i+\partial_i v)^2dx_1 dx_2+\sum\limits_{i,j=1}^2b_{ij}\int\limits_S (\partial_i\varphi_j+\partial_j\varphi_i)^2dx_1 dx_2
\end{eqnarray}
and
\begin{eqnarray}
\label{ener_mt}
 I(u,\varphi, v)=\sum\limits_{i=1}^2a_i\int\limits_S (\varphi_i+\partial_i v)^2dx_1 dx_2+\sum\limits_{i,j=1}^2b_{ij}\int\limits_S (\partial_i\varphi_j+\partial_j\varphi_i)^2dx_1 dx_2\nonumber\\+\sum\limits_{i,j=1}^2c_{ij}\int\limits_S (\partial_iu_j+\partial_ju_i+\frac{1}{2}\partial_iv_j\partial_jv_i)^2dx_1 dx_2,
\end{eqnarray}
where $v(x_1,x_2)$ is the (renormalized) vertical displacement and $u_i(x_1,x_2)$ are the (renormalized) in-plane displacements of the midsurface and $\varphi_i(x_1,x_2)$ are the angles of deflection of the orthogonal cross-sections of the plate. 
On $S\times(-h/2,h/2)$, this corresponds to a deformation of the form (we use the notation $x':=(x_1,x_2)$)
\begin{eqnarray}\label{eq:ansatz}y^h(x',x_3):=\left(\begin{array}{c} x'\\x_3   \end{array}\right)+\left(\begin{array}{c} h^{\alpha}u(x')\\h^{\beta}v(x')   \end{array}  \right)+x_3 \left(\begin{array}{c} h^{\gamma}\varphi(x')\\ 0\end{array}  \right).
\end{eqnarray}
The corresponding strain can be seen to be of the form
\begin{eqnarray*}
    (\nabla y^h)^T(\nabla y^h)= Id + \begin{pmatrix}
        &2 h^{\alpha} (\nabla' u)_{\text{sym}} + 2 x_3 h^{\gamma} (\nabla' \varphi)_{\text{sym}} +  h^{2\beta} \nabla' v \otimes \nabla' v & h^{\gamma} \varphi + h^{\beta} \nabla' v \\ &(h^{\gamma} \varphi + h^{\beta} \nabla' v)^T & 0
    \end{pmatrix} + \text{ h.o.t},
\end{eqnarray*}
where h.o.t. stands for higher order terms.
Hence, for an energy density $W$ that has the typical behavior $W(F) \sim |F^TF - \text{Id}|^2$ close to $SO(3)$ we find
\begin{align} \label{eq: scaling energy}
&\int_{S \times (-h/2,h/2)} W(\nabla y^h) \, dx \\ \sim & h \int_S \left( 2 \left| h^{\gamma} \varphi + h^{\beta} \nabla' v \right|^2 + \frac13 h^{2\gamma+2} \left| (\nabla' \varphi)_{\text{sym}} \right|^2 + 4 \left| h^{\alpha} (\nabla' u)_{\text{sym}} + \frac12 h^{2\beta}\nabla' v \otimes \nabla v\right|^2 \right) \, dx_1\,dx_2. \nonumber
\end{align}
Comparing with \eqref{ener_mtl} and \eqref{ener_mt} we notice the following:
\begin{itemize}
    \item it must hold $\gamma = \beta$ and $\alpha = 2\beta$;
    \item a renormalization of \eqref{eq: scaling energy} by $h^{2\gamma+3}$ will allow to obtain the term $\int_{S} |(\nabla' \varphi)_{\text{sym}}|^2 \, d\mathcal{L}^2$ in the limit $h\to 0$ but necessarily forces $\varphi = \nabla' v$, c.f., for example, \cite{FGM_Der};
    \item in order to derive all terms in \eqref{ener_mt} or \eqref{ener_mtl} from our ansatz the elastic energy density has to weight different entries of the elastic strain with different powers in $h$.
\end{itemize}
Taking into account the above considerations, we will derive  \eqref{ener_mtl} and \eqref{ener_mt} from an elastic energy that weights different entries of the elastic strain differently.
In our setting (see Section \ref{sec: notation} for the precise definitions) the derivation of \eqref{ener_mt} will be performed for $\beta = \gamma > 1$ whereas \eqref{ener_mtl} will be derived for $\gamma = \beta = 1$ and $\alpha = 2 \beta = 2\gamma +2 =4$. \\

Let us briefly discuss some related  literature.  Falach, Paroni, Podio-Guidugli and Tomasetti \cite{FPP}, \cite{PPT, PPT2} used an anisotropic (transversely isotropic) linear three-dimensional energy containing second-gradient terms. In fact, the idea was to consider for 
a given three-dimensional problem (the so-called “real problem”) an approximation whose variational limit coincides with the variational limit of the “real
problem”. Another idea is to consider anisotropic elasticity and to scale the elastic constants in  different ways, which  allows to avoid the inconsistency with the fact that for the Mindlin-Timoshenko-type models the shear modulus cannot be "too large". However, the parent models in \cite{FPP,PPT}  consider linearized strains and the Saint Venant–Kirchhoff model for the strain energy, which does not cover the case of nonlinear Lagrangian  strains and important models such as  neo-Hookean, Mooney-Rivlin, Ogden, Fung etc. (see \cite{Ogden} and references therein). 
Another approach starting from isotropic, linear elasticity with microrotations is presented in \cite{nhj}. For further discussion of the challenges and related literature we refer to the above references. \\\par

\par The rest of the paper is organized as follows: in Subsection \ref{sec: notation}  we set the main notation, state the parent problem and formulate the main assumptions on the energy density and external forces, we also present an example of a nonlinear energy functional satisfying these assumptions. In Subsection \ref{sec:mainresult} we formulate and discuss our main result. In Section \ref{sec:geometricrigidity}  rigidity estimates for the displacement gradient to the case of anisotropic elasticity are adjusted. In Section \ref{sec:dimred} we prove main results on $\Gamma$-convergence of the three-dimensional problem to the Reissner-Mindlin energy functional in case of the energy scaling of power greater than 4. Finally, Section \ref{sec:dimred2} contains a $\Gamma$-convergence result in case of  the scaling of  power 4.

\subsection{Notation and setting.} \label{sec: notation}
Throughout the text, we denote by $C$ generic constants that may change from expression to expression. \\
Let $S\subseteq\R^2$ be an open, bounded Lipschitz  domain, and consider for $h>0$ the (thin) domain 
\[\Omega_h=S\times \left( -\frac{h}{2}, \frac{h}{2}\right).\]
We will always assume that $h\in (0,1)$. For a deformation $w\in W^{1,2}(\Omega_h; \mathbb R^3)$, we consider the elastic energy
\begin{eqnarray}
\label{energy_func1}
\mathcal E_1^h(w):=\int\limits_{\Omega_h}W_1(\nabla'w'(z)^T \nabla'w'(z)+\nabla'w_3(z) \otimes \nabla'w_3(z))dz+h^{2}\int\limits_{\Omega_h}W_2(\nabla w(z)) dz,
\end{eqnarray}
where we use the notation $x'=(x_1, x_2)$, $w'=(w_1, w_2)$ and $\nabla'w=w,_1\otimes e_1+w,_2\otimes e_2$. Here and in the following, $w_{,i}$ denotes the $i$-th partial derivative $\partial_i w$ for $i=1,2,3$, and similarly for second order partial derivatives. For $M=(m_{ij})_{i,j=1\dots 3}\in\R^{3\times 3}$, we similarly set $M':=(m_{ij})_{i,j=1,2}\in \R^{2\times 2}$.
Now, consider the fixed domain $\Omega=S\times I$, where $I:=(-\frac12, \frac12)$ and introduce the change of variables $z(x)=(x_1,x_2,\frac{x_3}{h})$ and the rescaled deformation $y:\Omega\to \mathbb R^3$, $y(x)=w(z(x))$. It then holds  
$$
\nabla w=\left( \nabla' y, \frac{1}{h}y,_3\right) =:\nabla_h y
$$
and 
\begin{eqnarray}
\label{energy_func2}
\frac{1}{h}\mathcal E^h(y):=\int\limits_\Omega W_1(\nabla'y'(x)^T \nabla'y'(x)+\nabla'y_3(x) \otimes \nabla'y_3(x))dx+h^{2}\int\limits_{\Omega}W_2(\nabla_h y(x)) dx.
\end{eqnarray}
Slightly abusing notation, we identify functions $f:S\to\R^n$ with their trivial extensions $f:\Omega\to\R^n$ given by $f(x',x_3):=f(x')$. We sometimes write $u(x')$ or $u(x_3)$ instead of $u$ to point out on which components the function depends.\\

\textbf{Free energy densities.} 
We will always impose (without further mentioning) that the following assumptions on the free energy densities $W_i: \R^{(i+1)\times (i+1)} \to [0,\infty]$, $i=1,2$ hold. 
\begin{assumption}\label{ass:W}
    There exist constants $c_0, c_1, c_2>0$ such that 
\begin{enumerate}[label={(A\arabic*)}]
    \item\label{A1} $W_2(QF)=W_2(F)$ for all $Q\in SO(3)$ and all $F\in \R^{3\times 3}$,
    \item\label{A2} $W_1(Id)=\min W_1=0=\min W_2=W_2(Q)$ for all $Q\in SO(3)$,
    \item\label{A3} $W_2(F)\ge c_0 \,\dist^2 (F, SO(3))$ for all $F\in \R^{3\times 3}$,
    \item\label{A4}  $W_1(F'^TF'+(f_{31},f_{32})\otimes (f_{31},f_{32}))+c_1|(f_{31},f_{32})|^4\geq c_2\dist^2 (F', SO(2))$ for all $F=(f_{ij})\in\R^{3\times 3}$,  and
    \item\label{A5} $W_2$ is twice continuously differentiable in a neighborhood of $SO(3)$, and $W_1$ is twice differentiable in a neighbourhood of $\text{Id}$.
\end{enumerate}
\end{assumption}
Let us briefly discuss an example of free energy densities $W_1$ and $W_2$ satisfying Assumption \ref{ass:W}.
\begin{example}
We consider a special case of an orthotropic neo-Hookian strain energy (see e.g \cite{BSN, BB, SN}). Precisely, for $\beta>0$ we set 
$W_1:\R^{2\times 2}\to[0,\infty]$,
 \begin{equation}\label{ww1}W_1(C'):=\begin{cases}
    \beta \left[\tr\left(C'^2\right)-2\tr\left(C'\right)+2 \right]&\text{\qquad if\quad}C'\in\calA\\
+\infty&\text{\qquad otherwise,}
    \end{cases}
    \end{equation}
    where  
    \[\calA:=\left\{C'\in\R^{2\times2}:\ \exists F\in\R^{3\times 3}\text{ s.t. }\det F'>0\text{\ and }C'=F'^TF'+(f_{31},f_{32})\otimes (f_{31},f_{32})\right\}.\] 
In addition, for $\lambda,\mu>0$, we set $W_2:\R^{3\times 3}\to [0,\infty]$,
\[W_2(F):=\begin{cases}    
\mu(\tr (F^TF)-3)-\mu \ln \det (F^TF)+\lambda(\det (F^TF) -1)^2&\text{\qquad if\quad}\det F>0\\
+\infty&\text{\qquad otherwise.}
\end{cases}
\]
We note that the assumption $\det F>0$ is related to non-interpenetrability of matter, and the assumption $\det F'>0$ is typically satisfied in the case of infinitesimal planar strain.\\
Clearly, the functions $W_1$ and $W_2$ satisfy assumptions \ref{A1} and \ref{A5}. To see that $W_1$ satisfies \ref{A2} note that all $C'\in\calA$ are symmetric, and hence, denoting the eigenvalues of $C'$ by $\nu_1$ and $\nu_2$, we have 
\[W_1(C')=\beta\left[\nu_1^2+\nu_2^2-2(\nu_1+\nu_2)+2\right]=\beta \left[(\nu_1-1)^2+(\nu_2-1)^2\right]\geq 0\]
with equality if and only if $\nu_1=\nu_2=1$.\\
To see the other properties, we note that the matrix $F^TF$ is symmetric and positive definite and hence 
has three positive real eigenvalues that we denote by $\lambda_1^2$, $\lambda_2^2$ and $\lambda_3^2$ with $\lambda_i>0$ for $i=1,\dots,3$.
 Then 
 \begin{eqnarray*}
W_2(F^TF)&\ge& \mu(\lambda_1^2+\lambda_2^2+\lambda_3^2-3)-\mu (\ln \lambda_1^2+\ln \lambda_2^2+\ln \lambda_3^2)\ge
\mu\sum\limits_{i=1}^3(\lambda_i^2-\ln \lambda_i^2-1)\\
&\ge&
\mu\sum\limits_{i=1}^3(\lambda_i-1)^2=\mu|\sqrt{F^TF}-I|^2=\mu\, \dist^2(F, SO(3)),
\end{eqnarray*}
where we used that $\ln \lambda_i^2 = 2 \ln \lambda_i \leq 2 (\lambda_i - 1)$. This shows \ref{A3} and the first part of \ref{A2}.
Next, we consider $W_1$. For $F\in\R^{3\times 3}$, we have 
\begin{eqnarray}\label{eq:ex1}
&&W_1\left(F'^TF'+(f_{31},f_{32})\otimes (f_{31},f_{32})\right)\nonumber\\&\ge&  
\beta\left[
\left((f_{11}^2+f_{21}^2+f_{31}^2)^2+2(f_{11}f_{12}+f_{21}f_{22}+f_{31}f_{32})^2+(f_{12}^2+f_{22}^2+f_{32}^2)^2\right)\right.\nonumber\\
&&\left.-2\left(f_{11}^2+f_{21}^2+f_{31}^2+f_{12}^2+f_{22}^2+f_{32}^2\right)+2\right] \nonumber\\
&\geq& \beta\left[
(f_{11}^2+f_{21}^2-1)^2+2(f_{11}f_{12}+f_{21}f_{22})^2+(f_{12}^2+f_{22}^2-1)^2 \right.\nonumber\\&&\left.+f_{31}^4+f_{32}^4+2f_{31}^2f_{32}^2+2(f_{11}^2+f_{21}^2-1)f_{31}^2+2(f_{12}^2+f_{22}^2-1)f_{32}^2\right]\nonumber\\
&\geq&c_2 |F'^T F'-I|^2-c_1|(f_{31},f_{32})|^4.
\end{eqnarray}
\end{example}
Let us briefly explain how the previous example arises in modelling (see e.g \cite{BSN, BB, SN}). 
Orthotropic
materials are characterized by symmetry relations with respect to three orthogonal planes. The corresponding preferred directions are chosen as the intersections of these planes and are given by unit vectors
$a$, $b$ and $c$. If the vectors $a$, $b$ and $c$ are oriented along the coordinate axes, we can introduce the matrices 
$$M_1=\left(\begin{array}{ccc}
     1 &  0& 0\\
       0 &  0& 0\\
       0 &  0& 0\\
 \end{array}  \right),\qquad\qquad M_2=\left(\begin{array}{ccc}
     0&  0& 0\\
       0 &  1& 0\\
       0 &  0& 0\\
 \end{array}  \right)$$
 for structural tensors of the material, corresponding to the first two directions.
With the notation $F=\nabla y$, $F'=\nabla' y'$, ${\bf C}=F^TF$, 
the orthotropic part has the form
 \[W_1^{(M_1,M_2)}({\bf C})=\begin{cases}
    \beta \left[\sum\limits_{i=1}^2(\tr {\bf C}^2M_i-2\tr {\bf C}M_i +1)\right]&\text{\qquad if\quad}\det F'>0\\
+\infty&\text{\qquad otherwise.}
    \end{cases}\] 
 Since for $C'=F'^T F'+(f_{31},f_{32})\otimes (f_{31},f_{32})$, there holds
$$\tr C'^2=\tr C^2(M_1+M_2),\qquad \tr C'=\tr C(M_1+M_2),$$
we find that $W_1^{(M_1,M_2)}$ defined in this way coincides with \eqref{ww1}.\\[.2cm]

\textbf{Regularized functionals and constraints.}  We impose Dirichlet boundary conditions on admissible deformations, and following \cite{PPT},  we consider for technical reasons regularized versions of \eqref{energy_func2}. Precisely, we set $I^h:W^{1,2}(\Omega;\R^3)\to[0,\infty]$, 
\begin{eqnarray}
\label{func_i}
I^h(y):=\begin{cases}
   \int_\Omega W_1\left(\nabla'y'(x)^T \nabla'y'(x)+\nabla'y_3(x) \otimes \nabla'y_3(x)\right)dx+c_1\int\limits_{\Omega} |\nabla'y_3(x)|^4 dx+\\+h^2\int_{\Omega}W_2(\nabla_h y(x)) dx+\frac{1}{h^\epsilon}\int\limits_{\Omega} |y',_{33}(x)|^2 dx\text{\qquad if $y(x)=(x', hx_3)$ on $\partial S\times (0,1)$, and }\\
  +\infty\text{\qquad otherwise}
  \end{cases} 
\end{eqnarray}
with some $\epsilon > 0$. Here, the first expression is understood in the sense that it is $+\infty$ if $\nabla'y_3\not\in L^4$ or $y_{,33}'\not\in L^2$.
We note that we impose the specific boundary condition  on all of $\partial S$ in \eqref{func_i} only for the ease of notation. It can be easily relaxed to hold only on part of the boundary. We will point out during the proof explicitly which parts require the boundary condition and which ones also hold without it.

Note that the last term in \eqref{func_i} coincides with the second-gradient term in the linear three-dimensional functionals in \cite{PPT}.\\

\textbf{Forces.} To include forces $f^h\in L^2(S;\R^3)$, we follow \cite{FGM_Der} and assume that the total force and the total moment applied to the reference configuration is zero, i.e.
\begin{equation}
    \label{force_a1}
    \int\limits_\Omega f^h dx=0,\qquad\qquad\qquad     \int\limits_\Omega x\wedge f^h dx=0.
\end{equation}
We suppose that there is some $\alpha\in \R$ such that 
 \begin{equation}
   \label{force}
\frac{1}{h^\alpha}f^h\rightharpoonup f  \qquad \text{in} \quad L^2(S;\mathbb R^3) \text{\qquad with\qquad} f_1=f_2=0,
    \end{equation}
and denote the functional with forces by $J^h:W^{1,2}(\Omega;\R^3)\to(-\infty,\infty]$, 
\begin{equation}
\label{func_j}
J^h(y)=I^h(y) -\int\limits_\Omega f^h(x')\cdot y\, dx.
\end{equation}
\subsection{Main results and discussion.}\label{sec:mainresult}
Our main result is the $\Gamma$-limit of the sequence of functionals \eqref{func_i} as $h\to 0$. It turns out that the $\Gamma$-limit has the form of the linearized Reissner-Mindlin energy \eqref{ener_mt}, and coincides with the $\Gamma$-limit of \eqref{energy_func2}. Precisely, we have the following result.
  \begin{theorem}
      \label{th_conv1}
      \begin{enumerate} 
\item[(i)]{\em Compactness and lower bound.}  \\
 Suppose that $\sigma\geq 4$. Let $(y^h)_{h>0}\subseteq W^{1,2}(\Omega; \R^3)$ be such that 
 \begin{equation}
  \label{boun}
\limsup\limits_{h\to 0}\frac{1}{h^\sigma}I^h(y^h)<\infty.
  \end{equation}
  Then there exist  $Q_h:S\to SO(2)$, $u\in W^{1,2}(S;\R^2)$, $v\in W^{1,2}(S)$, and $\varphi\in W^{1,2}(S;\R^2)$ identified with the constant in $x_3$-direction function $\varphi\in W^{1,2}(\Omega;\R^2)$ such that
   \begin{align}
       &u^h\rightharpoonup  u\quad  \text{in} \;\; W^{1,2}(S;\mathbb R^{2}), \text{\quad where\quad}u^h(x') :=\frac{1}{h^{\sigma/2}}\int\limits_I\left(\left(\begin{array}{c} y_1^h(x',x_3)\\y_2^h(x',x_3)\end{array}\right)-x'\right) dx_3\label{u_conve},  \\
       & v^h\rightharpoonup  v\quad \text{in} \;\; W^{1,2}(S)\label{v_conve} ,\text{\quad where\quad} v^h(x'):=\frac{1}{h^{\sigma/2-1}} \int\limits_I y_3^h(x',x_3)\,dx_3, \\
      & \varphi^h\rightharpoonup \varphi \quad  \text{in} \;\,   L^{2}(\Omega;\mathbb R^2),\text{\quad where\quad} \varphi^h(x',x_3):=\frac{1}{h^{\sigma/2}}\left(Q_h^T\partial_3 y'(x',x_3)\right),\label{phi_conve1}  \\
      & \frac{1}{h^{\sigma/2}}Q_h^T\left(\nabla' (y^h)'-\int\limits_I \nabla' (y^h)'dx_3\right)\rightharpoonup x_3\nabla'\varphi \; \text{in} \;   L^{2}(\Omega;\mathbb R^{2})\label{phi_conve2}  .
   \end{align}
 Moreover, 
if $\sigma>4$
set
\begin{equation*}      \Q_3^i(A):=\frac{\partial^2W_i(A)}{\partial A^2}(Id)(A,A), \quad i=1,2,\text{\quad and\quad}
\Q_2^2(\tilde G):=\min_{c\in \mathbb R}\\
\Q_3^2\left(\tilde G+c e_3\otimes e_3\right).
\end{equation*}
Then 
we have 
  \begin{eqnarray}
      \label{inflim1}
     \liminf\limits_{h\to 0}\frac{1}{h^\sigma}I^h(y^h)\ge
\frac{1}{2}\int\limits_S  \Q_2^2(\tilde G)\,dx'+
      \frac{1}{2}\int\limits_S  \Q_3^1(2\sym \nabla'u)\,dx'+\frac{1}{6}\int\limits_S \Q_3^1(\sym \nabla'\varphi)\,dx'
 \end{eqnarray}    
 with
      \begin{eqnarray}
          \label{eq:tildeG}
 \tilde G(x'):=\left( \begin{array}{ccc}
     0&0&\varphi_1(x')+\partial_1 v(x')  \\
     0&0&\varphi_2(x') +\partial_2 v(x') \\
     \varphi_1(x')+\partial_1 v(x')&\varphi_2(x') +\partial_2 v(x')&0
\end{array}\right).
  \end{eqnarray}
 \item[(ii)]{\em Optimality of lower bound in the case $\sigma>4$.}\\
 If  $u\in W^{1,2}(S;\R^2)$, $v\in W^{1,2}(S)$, and $\varphi\in W^{1,2}(S;\R^2)$, then there exists a recovery sequence $(\hat y^h)_{h>0}$ such that \eqref{u_conve}--\eqref{phi_conve2} hold and 
 \begin{eqnarray}
     \label{opt}
      \liminf\limits_{h\to 0}\frac{1}{h^\sigma}I^h(\hat y^h)= 
\frac{1}{2}\int\limits_S  \Q_2^2(\tilde G)dx'+
       \frac{1}{2}\int\limits_S \Q_3^1(2\sym \nabla'u)dx'+\frac{1}{6}\int\limits_S  \Q_3^1(\sym \nabla'\varphi(x'))dx'.
 \end{eqnarray}    
 \end{enumerate}
  \end{theorem} 
As a corollary of Theorem \ref{th_conv1} one can infer convergence of minimizers to minimizers of a linear Reissner-Mindlin energy. Precisely, we have the following result.
  \begin{theorem}
\label{th_j}
Suppose that $\sigma>4$. Consider external forces $(f^h)_{h>0}$ such that \eqref{force_a1} and \eqref{force} hold with $\sigma=2\alpha-2$. Let $(y^h)_{h>0}$  be a $\sigma$-minimizing sequence for $(J^h)_{h>0}$, i.e.,
$$\limsup\limits_{h\to 0} \frac{1}{h^\sigma}\left(J^h(y^h)-\inf J^h\right)=0.$$
Then there exist  $Q_h:S\to SO(2)$, $u\in W^{1,2}(S;\R^2)$, $v\in W^{1,2}(S)$, and $\varphi\in W^{1,2}(S;\R^2)$  such that 
$$\nabla_hy^h\to Id \quad \text{in}\quad L^2(\Omega;\mathbb R^3),$$
\eqref{u_conve} holds with $u=0$, and \eqref{v_conve}--\eqref{phi_conve1} hold up to a subsequence. Moreover, $0\ge \inf J^h\ge -Ch^\sigma$ and the limit function $(v,\varphi)$ 
minimizes the functional
\begin{eqnarray*}
J(v,\varphi)=\frac{1}{2}\int\limits_S  \Q_2^2(\tilde G)\,dx'+\frac{1}{6}\int\limits_S \Q_3^1(\sym \nabla'\varphi(x'))\,dx'-\int\limits_S f_3(x') v(x') \,dx'
\end{eqnarray*}
with $\tilde{G}$ defined in \eqref{eq:tildeG}. 
Finally,  $\lim_{h\to0}\frac{1}{h^\sigma}J^h(y^h)=\min J.$  
\end{theorem}
Let us briefly discuss the main difficulties in our analysis. Starting point is the seminal work in \cite{FGM_Der} and \cite{FGM_Rig}on dimension reduction in nonlinear elasticity. However, the scaling of the energy functional considered there yields a higher rigidity than we expect in the Reissner-Mindlin model. In particular, one obtains relations between the angles of rotation of the cross-sections of the plate and the transversal displacement of the mid surface $(\phi_1,\phi_2)=-\nabla v$. The Reissner-Mindlin model on the other hand takes into account transverse shear effects, which contradict this relation. \\
To overcome this, we rescale parts of the energy differently. More precisely, if we denote the strain matrix by $S=(s_{ij})_{i,j=1,2,3}$, we rescale the part depending on the upper part $S'=(s_{ij})_{i,j=1,2}$ by $h^\sigma$ while the rest of the energy is rescaled by $h^{\sigma-2}$. While this rescaling allows us to avoid unwanted rigidity, it also leads to a lack of compactness. In case $\sigma>4$, one can overcome this difficulty by adding regularizing terms (see \eqref{func_i}), which vanish in the limit. Here we use 
Dirichlet-type boundary conditions, which allow us to use Korn-type inequalities (see e.g. Proposition 1 \cite{FGM_Der}). \\However, for the case $\sigma=4$, we still experience the lack of compactness to perform the limit transition in the nonlinear term $\nabla'y_3(x) \otimes \nabla'y_3(x)$ in \eqref{func_i}. The way to overcome this difficulty by adding a second gradient term is discussed in Section 4 but  the limit functional differs from the classical Reissner-Mindlin energy.\\
The main steps in the proof are the following:
\begin{enumerate}
    \item[(a)] {\it Scaled  rigidity estimates in thin domain.} In a thin domain $\Omega_h=S\times(-h/2,h/2)$ we use different rigidity estimates. We derive approximations of the scaled deformation gradient $\nabla_h y$ by a rotation $R_h(x')$ in $SO(3)$ depending on $x'$, the averaged with respect to the third component 2D gradient $\int\limits_I\nabla' y'dx_3$ by a rotation $Q_h(x')$ in $SO(2)$ depending on $x'$ (see Theorem \ref{th_rig1}), and $\nabla'y'$ by a rotation $T_h(x)$ in $SO(2)$ depending on $x$ (see Theorem \ref{th_rig2}). In the first case the approximation rate is $h^{\sigma-2}$ and in the second and the third cases  $h^{\sigma}$. We also construct a rotation $L_h\in SO(3)$ associated with $Q_h$
    $$L_h=\left(\begin{array}{ccc}Q_h& &0\\  & &0\\0&0&1    \end{array}\right).$$
    \item[(b)] {\it Scaling of deformations.}
    Next we obtain the rates of convergence of rotations to the identity in $L^2$ and in any $L^p$ for $2<p<\infty$ and  normalize and scale in-plane and out-of-plane deformations according to these estimates. We note that the scaling of the out-of-plane components does not coincide with that introduced in \cite{FGM_Der}. 
    \item[(c)] {\it Convergence of scaled deformations.}
    We prove weak convergence of scaled and averaged with respect to the third component deformations. Using the boundary conditions we show that the deformations converge to the same functions even without normalization.
    \item[(d)] {\it Identification of the limiting strain.}
    We then estimate in $L^2$ three approximate strain components and find that their weak limits in $L^2$ up to a subsequence have the following structures:
    $$G_h:=\frac{L_h^T\nabla_h y^h -Id}{h^{\frac{\sigma}{2}-1}}\rightharpoonup \left( \begin{array}{ccc}
     0&0&\varphi_1(x')  \\
     0&0&\varphi_2(x')  \\
     \partial_1 v(x')&\partial_2 v(x')&G_{33}(x).
\end{array}\right)$$
$$\sym  F_h:=\sym \frac{Q_h^T\int\limits_I\nabla' (y^h)' dx_3-Id}{h^{\frac{\sigma}{2}}}\rightharpoonup  \sym \nabla'u(x')$$
$$K_h:=\frac{Q_h^T(\nabla' (y^h)'-\int\limits_I\nabla' (y^h)' dx_3)}{h^\frac{\sigma}{2}}\rightharpoonup x_3\nabla'\varphi(x').$$
    \item[(e)]{\it $\Gamma$-convergence. }
    Using the Taylor's expansion and this relations we obtain the lower bound in Theorem \ref{th_conv1}. For the proof of the optimality of the lower bound we construct a recovering sequence \eqref{recov}.
\end{enumerate}

\section{Geometric rigidity}\label{sec:geometricrigidity}
To prove the results on the approximation of the deformation gradient by rotations we make use of the celebrated Friesecke-James-M\"uller rigidity theorem: 
\begin{theorem}[\cite{FJM}]
\label{th_rig_gen}
Let $U$ be a bounded Lipschitz domain in $\mathbb R^n$,  $n\ge 2$ . There exists a constant $C(U)$ such that for each $v\in W^{1,2}(U, \mathbb R^n)$ there is an associated rotation $R\in SO(n)$ satisfying 
$$
\|\nabla v-R\|_{L^2(U)}\le C(U)\|\dist(\nabla v, SO(n)\|_{L^2(U)}.
$$
\end{theorem}
\begin{remark}
    The constant $C(U)$ in Theorem \ref{th_rig_gen} can be chosen uniformly for a family of domains which are Bilipschitz equivalent with controlled Lipschitz constant (see also \cite[Section 5]{conti-zwicknagl}). The constant $C(U)$ is invariant under dilations.
\end{remark}
Building on this estimate, we show the following result.
\begin{theorem}
\label{th_rig1}
Let $S\subseteq \mathbb R^2$ be a Lipschitz domain and $\Omega :=S\times (-\frac{1}{2}, \frac{1}{2})$. Let $y\in W^{1,2}(\Omega, \mathbb R^3)$
and set
\begin{equation}
\label{e}
E_1:=\int\limits_\Omega \dist^2 (\nabla_h y, SO(3)) dx\le Ch^2,\text{\ and\ }\; E_2:=\int\limits_\Omega \dist^2 (\nabla' y', SO(2)) dx\le Ch^2.
\end{equation}
 Then there exist 
maps  $ R_h:S\to SO(3)$, $\tilde R_h:S\to \mathbb R^{3\times 3}$,  
$ Q_h:S\to SO(2)$, and $\tilde Q_h:S\to \mathbb R^{2 \times 2}$ such that $|\tilde R_h|\le C$, $|\tilde Q_h|\le C$ and $\tilde R_h\in W^{1,2}(S,\mathbb R^{3\times 3})$, $\tilde Q_h\in W^{1,2}(S,\mathbb R^{2\times 2})$. Moreover,
\begin{align}
\|\nabla_h y(x',x_3)-R_h(x')\|_{L^2(\Omega)}^2\le CE_1, \quad\|\tilde R_h(x')-R_h(x')\|_{L^2(S)}^2\le CE_1, \label{r1}\\
\|\nabla'\tilde R_h(x')\|_{L^2(S)}^2\le \frac{C}{h^2}E_1, \quad\|\tilde R_h(x')-R_h(x')\|_{L^\infty(S)}^2\le \frac{C}{h^2}E_1, \label{r2}\\
\|\int\limits_I\nabla' y'(x)dx_3-Q_h(x')\|_{L^2(S)}^2\le CE_2, \;\|\tilde Q_h(x')-Q_h(x')\|_{L^2(S)}^2\le CE_2, \label{q1}\\
\|\nabla'\tilde Q_h(x')\|_{L^2(S)}^2\le \frac{C}{h^2}E_2, \quad\text{\ and\quad}\|\tilde Q_h(x')-Q_h(x')\|_{L^\infty(S)}^2\le \frac{C}{h^2}E_2. \label{q2}
\end{align}
\end{theorem}
\begin{proof}
Estimates \eqref{r1} and \eqref{r2} can be derived as in \cite[Theorem 6]{FGM_Der}, so we omit the proof here. \par
Let $U$ be an open subset in $\mathbb R^2$ and let $K\subset U$ be compact and such that $\dist_\infty (K, \partial U)>3h$, where $\dist_\infty$ is the distance with respect to the norm  $\|(x_1, x_2)\|_\infty=\text{max} \{|x_1|, |x_2|\}$. For each point $x'\in K$ we consider the square 
$$  S_{x',h}:= x'+(0,h)^2$$
with lower left corner  $x'$.
Let $\psi\in C_c^\infty((0, 1)^2)$ be a standard mollifier, i.e. $\psi\ge 0$ and $\int_{\R^n} \psi =1$, and set $\psi_h(\cdot)=h^{-2}\psi(\cdot/h)$. We set 
\begin{equation}
\label{notatf}
F(x', x_3):=\nabla' y'(x',x_3)  \quad \text{and}  \quad \bar F(x'):=\int\limits_I  F(x', x_3) dx_3
\end{equation}
and consider the map
$$\hat Q_{h}(x'):=(\psi_h * \bar F)(x')= \int\limits_{S_{x',h}}h^{-2} \psi\left( \frac{x'-z'}{h}  \right) \int\limits_I F(z',z_3)dz_3 dz'.$$
Applying Theorem \ref{th_rig_gen} to $S_{x',h}$ we obtain that for any fixed $z_3\in I$ there exists a rotation $Q_{x',z_3, h}$  such that
\begin{equation}
\label{s_rig}
\int\limits_{S_{x',h}}|F(z', z_3)- Q_{x',z_3, h}|^2 dz'\le C \int\limits_{S_{x',h}} \dist^2(F(z',z_3), SO(2)) dz'.
\end{equation}
Hence, by Hölder's inequality and Fubini, we also obtain 
\begin{eqnarray}
\label{s_i_rig}
&&\int\limits_{S_{x',h}}\left |\int\limits_I F(z', z_3)d z_3- \int\limits_I Q_{x',z_3, h}d z_3 \right|^2 dz'\nonumber\\&\le& \int\limits_I \int\limits_{S_{x',h}}|F(z', z_3)- Q_{x',z_3, h}|^2 dz' d z_3\le C \int\limits_{S_{x',h}\times I} \dist^2(F(z), SO(2)) dz.
\end{eqnarray}
Consequently, using Hölder's inequality and \eqref{s_i_rig} we have
\begin{eqnarray}
\label{hatr_rot}
&&\left|\hat Q_h(x')-\int\limits_I Q_{x',z_3, h}d z_3 \right|^2 \nonumber\\&\le& \left|\int\limits_{S_{x',h}} h^{-2}\psi\left( \frac{x'-z'}{h}  \right)\int\limits_I F(z', z_3) dz_3 dz'-\int\limits_{S_{x',h}} h^{-2}\psi\left( \frac{x'-z'}{h}  \right)\int\limits_I Q_{x',z_3,h} dz_3 dz' \right|^2\nonumber\\
&\le& \frac{C}{h^2} \int\limits_{S_{x',h}} \left |\int\limits_I F(z', z_3)d z_3- \int\limits_I Q_{x',z_3, h}d z_3 \right|^2 dz'
\nonumber\\
&\le& \frac{C}{h^2} \int\limits_{S_{x',h}\times I} \dist^2(F(z), SO(2)) dz.
\end{eqnarray}
Since $\int \nabla \psi_h=0$, for any point $\tilde x'\in S_{x',h}$ we have from Hölder's inequality and \eqref{s_i_rig}
\begin{eqnarray}   
\label{hatr_nabla1}
\left|\nabla \hat Q_h(\tilde x')\right|^2 &=& 
\left|\int\limits_{S_{\tilde x',h}} h^{-3 }\nabla \psi\left( \frac{\tilde x'-z'}{h}  \right)\int\limits_I F(z', z_3) dz_3 dz'-\int\limits_{S_{\tilde x',h}} h^{-3 }\nabla \psi\left( \frac{\tilde x'-z'}{h}  \right)\int\limits_I Q_{x',z_3,h} dz_3 dz'  \right|^2\nonumber\\
&\le& \frac{C}{h^4} \int\limits_{S_{x',2h}} \left |\int\limits_I F(z', z_3)d z_3- \int\limits_I Q_{x',z_3, h}d z_3 \right|^2 dz'
\nonumber\\
&\le& \frac{C}{h^4} \int\limits_{S_{x',2h}\times I} \dist^2(F(z), SO(2)) dz, 
\end{eqnarray}
and integrating \eqref{hatr_nabla1} over $S_{x',h}$ we find
\begin{equation}
\label{hatr_nabla2}
\int\limits_{S_{x',h}}\left|\nabla \hat Q_h(z')\right|^2 dz'\le  \frac{C}{h^2} \int\limits_{S_{x',2h}\times I} \dist^2(F(z), SO(2)) dz.
\end{equation}
For any point $\tilde x'\in S_{x',h}$ we have again using \eqref{s_i_rig}
\begin{eqnarray}
\label{hatr_hatr}
\left|\hat Q_h(x')- \hat Q_h(\tilde x') \right|^2 &\le& \left|\int\limits_{S_{x',h}} h^{-2}\psi\left( \frac{x'-z'}{h}  \right)\int\limits_I F(z', z_3) dz_3 dz'-\int\limits_{S_{\tilde x',h}} h^{-2}\psi\left( \frac{\tilde x'-z'}{h}  \right)\int\limits_I F(z', z_3) dz_3 dz'  \right|^2\nonumber\\&\le& 
C\left|\int\limits_{S_{x',h}} h^{-2}\psi\left( \frac{x'-z'}{h}  \right)\int\limits_I F(z', z_3) dz_3 dz'-\int\limits_{S_{ x',h}} h^{-2}\psi\left( \frac{ x'-z'}{h}  \right)\int\limits_I Q_{x',z_3, h} dz_3 dz'  \right|^2\nonumber\\&+&
C\left|\int\limits_{S_{\tilde x',h}} h^{-2}\psi\left( \frac{\tilde x'-z'}{h}  \right)\int\limits_I F(z', z_3) dz_3 dz'-\int\limits_{S_{ \tilde x',h}} h^{-2}\psi\left( \frac{ \tilde x'-z'}{h}  \right)\int\limits_I Q_{x',z_3, h} dz_3 dz'  \right|^2\nonumber\\
&\le&\frac{C}{h^2} \int\limits_{S_{x',2h}} \left |\int\limits_I F(z', z_3)d z_3- \int\limits_I Q_{x',z_3, h}d z_3 \right|^2 dz'\nonumber\\
&\le& \frac{C}{h^2} \int\limits_{S_{x',2h}\times I} \dist^2(F(z), SO(2)) dz.
\end{eqnarray}
Combining   \eqref{s_i_rig}, \eqref{hatr_rot}, and \eqref{hatr_hatr} we obtain
\begin{eqnarray}
\label{hatr_grad}
&&\int\limits_{S_{x',h}}\left|\hat Q_h(z')- \int\limits_I F(z', z_3) d z_3 \right|^2dz' \nonumber\\
&\le&C\left( \int\limits_{S_{x',h}}\left|\hat Q_h(z')- \hat Q_h(x')\right|^2dz'+
\int\limits_{S_{x',h}}\left|\hat Q_h(x')- \int\limits_I Q_{x',z_3,h} d z_3 \right|^2dz'+\int\limits_{S_{x',h}}\left| \int\limits_I Q_{x',z_3,h} d z_3 - \int\limits_I F(z', z_3) d z_3 \right|^2dz'\right)\nonumber\\
&\le& C \int\limits_{S_{x',h}\times I} \dist^2(F(z), SO(2)) dz.
\end{eqnarray}
Finally, setting $g(\zeta):=\dist (\hat Q_h(x'+h\zeta), SO(2))$ and taking into account \eqref{hatr_nabla1},  \eqref{hatr_grad}
we arrive at
\begin{eqnarray*}
&&\int_{(0,1)^2}|g|^2 d\zeta+\sup\limits_{(0,1)^2}|\nabla' g|^2
=\int\limits_{S_{x',h}} \dist^2(\hat Q_h(z'), SO(2)) dz'+\sup\limits_{S_{x',h}}|\nabla' \dist(\hat Q_h(z'), SO(2)) |^2\\
&\le& C\Big(\int\limits_{S_{x',h}} |\hat Q_h(z')- \int\limits_I 
 F(z',z_3)dz_3|^2dz'+\int\limits_{S_{x',h}} \dist^2(\int\limits_I 
 F(z',z_3)dz_3, SO(2)) dz'+\sup\limits_{S_{x',h}}|\nabla' \hat Q_h(z')|^2\Big)\\
&\le&\frac{C}{h^2} \int\limits_{S_{x',2h}\times I} \dist^2(F(z), SO(2)) dz.
\end{eqnarray*}
Consequently,
\begin{eqnarray}
\label{hatr_bound}
\dist^2 (\hat Q_h(\tilde x'), SO(2))\le \sup\limits_{(0,1)^2}|g|^2\le C(\int_{(0,1)^2}|g|^2 d\zeta+\sup\limits_{(0,1)^2}|\nabla' g|^2)\nonumber\\
\le\frac{C}{h^2} \int\limits_{S_{x',2h}\times I} \dist^2(F(z), SO(2)) dz,
\end{eqnarray}
and integrating \eqref{hatr_bound}  we obtain
\begin{eqnarray}
\label{hatr_bound2}
\int\limits_{S_{x',h}}\dist^2 (\hat Q_h(\tilde x'), SO(2)) d \tilde x'
\le\frac{C}{h^2} \int\limits_{S_{x',h}}\int\limits_{S_{x',2h}\times I} \dist^2(F(z), SO(2)) dz d\tilde x'\nonumber\\
\le
C \int\limits_{S_{x',h}\times I} \dist^2(F(z), SO(2)) dz.
\end{eqnarray}
Now we consider a lattice of squares of size $h$ in  $\mathbb R^2$ and sum \eqref{hatr_nabla2}, \eqref{hatr_grad} over all squares which intersect $K$. This yields
\begin{eqnarray}
\label{hatr_k}
\int\limits_{K}\left(|\hat Q_h(x')-\int\limits_I 
 F(x',x_3)dx_3|^2+h^2|\nabla' \hat Q_h(x')|^2\right)dx'\le C \int\limits_{U\times I} \dist^2(F(x), SO(2)) dx.
\end{eqnarray}
Following \cite{FGM_Der} we consider $S$ to be locally the epigraph of a Lipschitz function, then arguing as in \cite{FGM_Der} and applying the above estimates we get
\begin{eqnarray}
\label{hatr_k_s}
\int\limits_{K\cap S}\left(|\hat Q_h(x')-\int\limits_I 
 F(x',x_3)dx_3|^2+h^2|\nabla' \hat Q_h(x')|^2\right)dx'\le C \int\limits_{U\times I} \dist^2(F(x), SO(2)) dx.
\end{eqnarray}
Since $S$ is Lipschitz, its closure $\bar S$ can be covered by a finite number of open cubes $U_0,..., U_l$, where $\bar U_0\subset S$. Denote $\hat Q_i$ the maps constructed in previous steps for $U_i$ and consider a partition of unity corresponding to the cover $\{U_i\}$
$$ \eta_i\in C_0^\infty(U_i), \;\qquad\eta_i\ge 0, \;\qquad\sum\limits_{i=0}^l\eta_i(x')=1,\,\forall x'\in S.$$
If we set
\begin{equation}
\label{tilde_r}
\tilde Q_h=\sum\limits_{i=0}^l\eta_i \hat Q_i
\end{equation}
then 
\begin{align*}
\tilde Q_h-\int\limits_I F dx_3=\sum\limits_{i=0}^l\eta_i (\hat Q_i-\int\limits_I F dx_3) \text{\qquad and\qquad}
\nabla'\tilde Q_h=\sum\limits_{i=0}^l\eta_i \nabla'\hat Q_i+\sum\limits_{i=0}^l\nabla'\eta_i (\hat Q_i-\int\limits_I F dx_3).
\end{align*}
Consequently, applying \eqref{hatr_k} and \eqref{hatr_k_s} with $K_i=\text{supp} \,\eta_i$ we obtain
\begin{equation}
\label{tilder}
\int\limits_{S}\left(|\tilde Q_h(x')-\int\limits_I 
 F(x',x_3)dx_3|^2+h^2|\nabla' \tilde Q_h(x')|^2\right)dx'\le C E_2.
\end{equation}
Note that this implies in particular that $\|\nabla'\tilde{Q}_h(x')\|_{L^2(S)}^2\leq \frac{C}{h^2}E_2$, which is the first part of assertion of \eqref{q1}.
From \eqref{hatr_bound}  we obtain
\begin{eqnarray}
\label{tilder_bound}
\sup\limits_S\dist^2 (\tilde Q_h(x'), SO(2)) \le\frac{C}{h^2} \sup\limits_{x'\in S}\int\limits_{(B_{x', c_0h}\cap S)\times I} \dist^2(\nabla' y(z), SO(2)) dz\le \frac{C}{h^2} E_2\le C
\end{eqnarray}
where $c_0$ depends only on $S$. Note that this estimate implies in particular
\begin{eqnarray}
\label{eq:tildeQbounded}
\|\tilde{Q}_h\|_{L^\infty(S)}\leq C.
\end{eqnarray}
It follows from \eqref{hatr_bound2} that
\begin{equation}
\label{tilder_bound1}
\int\limits_{S}\dist^2 (\tilde Q_h(\tilde x'), SO(2)) d \tilde x'
\le CE_2.
\end{equation}
 Now the rotation $Q_h$ can be obtained as projection of $\tilde Q_h$ onto $SO(2)$. More precisely, since $SO(2)$ is a smooth manifold, there exists a tubular neighbourhood $\mathcal U$ of $SO(2)$
such that the projection $\pi:\mathcal U\to SO(2)$ is smooth. In particular, there exists some $\delta>0$ such that for all $M\not\in\mathcal{U}$, there holds $\dist(M,SO(2))>\delta$. Let
\begin{eqnarray}\label{eq:tildepi}
\tilde{\pi}:\R^{2\times 2}\to SO(2),\qquad \tilde \pi(M):=\begin{cases}
    \pi(M)&\qquad   M\in \mathcal U,
    \\  Id&\qquad M\notin \mathcal U,
\end{cases}
\end{eqnarray}
and set $Q_h:=\tilde{\pi}\circ Q_h$.
If $\tilde Q_h(x')\in\mathcal U$ then $|Q_h(x')-\tilde Q_h(x')|= \dist(\tilde Q_h, SO(2))$. If $\tilde Q_h(x')\notin\mathcal U$ then  $\dist(\tilde Q_h(x'), SO(2))>\delta$, and with \eqref{eq:tildeQbounded} we obtain $|Q_h(x')-\tilde Q_h(x')|\le C \dist(\tilde Q_h, SO(2))$. Consequently, with \eqref{tilder_bound} we obtain
\begin{equation}
\label{q_bound}
\|Q_h(x')-\tilde Q_h(x')\|_{L^\infty(S)}^2= \sup\limits_S|Q_h(x')-\tilde Q_h(x')|^2\le C \sup\limits_S\dist^2(\tilde Q_h, SO(2))\le \frac{C}{h^2} E_2,
\end{equation}
which shows the second assertion of \eqref{q2}. Similarly, with \eqref{tilder_bound1} we deduce
\begin{equation}
\label{q_bound1}
\int\limits_S|Q_h(x')-\tilde Q_h(x')|^2 dx'\le C \int\limits_S \dist^2(\tilde Q_h, SO(2)) dx'\le C E_2,
\end{equation}
which is the second assertion of \eqref{q1}. Finally, combining \eqref{tilder} and \eqref{q_bound1} we obtain the first assertion of \eqref{q1}. This concludes the proof of Theorem \ref{th_rig1}.
\end{proof}

\begin{remark}
\label{rem1}
Let $\mathcal{U}_1$ and $\mathcal{U}_2$ be neighborhoods of $SO(3)$ and $SO(2)$ such that the respective projections onto $SO(3)$ and $SO(2)$ are smooth. If we assume that $E_i\le Ch^{2+\varepsilon}$ for some $\varepsilon>0$ then (for $h$ small enough) we obtain from \eqref{r2} and \eqref{q2} that $\tilde R_h(x')\in \mathcal U_1$  and $\tilde Q(x')\in \mathcal U_2$. It follows that the the maps $R_h:S\to SO(3)$ and $Q_h:S\to SO(2)$ which are obtained by projections onto $SO(3)$ and $SO(2)$, respectively, are actually also in $W^{1, 2}$, and it holds
\begin{align}
\|\nabla' R_h(x')\|_{L^2(S)}^2\le \|\nabla' \tilde R_h(x')\|_{L^2(S)}^2\le \frac{C}{h^2}E_1,\label{r3}\\
\|\nabla' Q_h(x')\|_{L^2(S)}^2\le \|\nabla' \tilde Q_h(x')\|_{L^2(S)}^2\le \frac{C}{h^2}E_2. \label{q3}
\end{align}
Thus in this case the estimates from Theorem \ref{th_rig1} (more precisely the first assertions of \eqref{r2} and \eqref{q2}) hold for $R_h$ and $Q_h$ directly.
\end{remark}
\begin{corollary}
\label{cor1}
Under the assumptions of Theorem \ref{th_rig1} there exists a constant rotation $P_h'\in SO(2)$ such that 
\begin{align}
\|\int\limits_I\nabla' y'(x)dx_3-P_h'\|_{L^2(S)}^2\le \frac{C}{h^2} E_2, \text{\qquad and}\label{q5}\\
\quad\|Q_h(x')-P_h'\|_{L^p(S)}^2\le \frac{C}{h^2}E_2,\,\, \;1\le p<\infty .\label{q6}
\end{align}
\end{corollary}
\begin{proof}
We use the notation from Theorem \ref{th_rig1} and from \eqref{eq:tildepi} and define
$$P_h'=\tilde\pi\left(\frac{1}{|S|}\int_S\tilde Q_h(x') dx'\right). $$
It follows from Poincar\'{e}'s inequality that
\begin{align*}
    &\|\tilde Q_h(x')-P_h'\|_{L^2(S)}^2  \le C\left(\left\|\tilde\pi(\tilde Q_h(x'))-\tilde\pi\left(\frac{1}{|S|}\int_S\tilde Q_h(x') dx'\right)\right\|_{L^2(S)}^2 + \left\|\tilde Q_h(x')-\tilde\pi(\tilde Q_h(x'))\right\|_{L^2(S)}^2\right)\\ \le &C\left(\|\nabla' \tilde Q_h(x')\|_{L^2(S)}^2+ \int_S \dist^2(\tilde Q_h(x'),SO(2)) dx'  \right)\le C\frac{E_2}{h^2}.
\end{align*}
Hence by the Sobolev embedding, we have for $1 \leq p < \infty$
\begin{eqnarray*}\|\tilde Q_h(x')-P_h'\|_{L^p(S)}^2\le C\|\tilde Q_h(x')-P_h'\|_{W^{1,2}(S)}^2  \le C\left(\|\tilde Q_h(x')-P_h'\|_{L^2(S)}^2+\|\nabla'\tilde Q_h(x')\|_{L^2(S)}^2\right)\le \frac{C}{h^2}E_2, \end{eqnarray*}
where in the last step we used the first estimate in \eqref{q2}.
Then \eqref{q5} and \eqref{q6} follow with \eqref{q1}, \eqref{q2}.
\end{proof}
Now we define $R_h'$ as $2\times 2$ submatrix of the rotation matrix $R_h$,
and the rotation $P_h\in SO(3)$ via
\begin{equation}
\label{p_h}
R_h=\left(\begin{array}{ccc}R_h'& &0\\  & &0\\0 &0 & 1   \end{array}\right),\qquad\text{and\qquad}
P_h:=\left(\begin{array}{ccc}P_h'& &0\\  & &0\\0&0&1    \end{array}\right).
\end{equation}
\begin{corollary}
\label{cor2}
Under the assumptions of Theorem \ref{th_rig1} the constant rotation $P_h'\in SO(2)$ from Corollary \ref{cor1} satisfies
\begin{equation}
\label{r_p_prime}
\quad\|R_h'(x')-P_h'\|_{L^2(S)}^2\le C \left(\frac{E_2}{h^2}+E_1\right)
\end{equation}
and 
\begin{equation}
\label{r_p}
 \|R_{33}^h-1\|_{L^2(S)}^2\le C\left(\frac{E_2}{h^2}+E_1 \right).
\end{equation}
\end{corollary}
\begin{proof}
We obtain from \eqref{r1} and H\"older's inequality that
\[\|R_h'(x')-\int\limits_I\nabla' y'(x)dx_3\|_{L^2(S)}^2\leq \int_S\int_I\left|R_h'(x')-\nabla' y'(x)\right|^2\,dx_3dx',\]
and hence with \eqref{q5}, we deduce
\begin{align*}
\|R_h'(x')-P_h'\|_{L^2(S)}^2\le & C\left(\|R_h'(x')-\int\limits_I\nabla' y'(x)dx_3\|_{L^2(S)}^2 + \|\int\limits_I\nabla' y'(x)dx_3-P_h'\|_{L^2(S)}^2\right) \\ \le &C \left(\frac{E_2}{h^2}+E_1\right),
\end{align*}
which completes the proof of \eqref{r_p_prime}.\\
We now turn to the proof of \eqref{r_p}.
We denote the entries of the matrices $R_h$ and $P_h'$ by $R_{ij}^h$ and $P^h_{ij}$, respectively. 
Since $R_h$ and $P_h'$ are rotations it holds for $j=1,2$ that 
\begin{align*}
|R_{j1}^h|^2+|R_{j2}^h|^2+|R_{j3}^h|^2=1=
|P_{j1}^h|^2+|P_{j2}^h|^2.
\end{align*}
Consequently, for $j=1,2$
\begin{align*}
|R_{j3}^h|^2 & =|P_{j1}^h|^2-|R_{j1}^h|^2+|P_{j2}^h|^2-|R_{j2}^h|^2\\ &= 
(|P_{j1}^h|+|R_{j1}^h|)(|P_{j1}^h|-|R_{j1}^h|)+(|P_{j2}^h|+|R_{j2}^h|)(|P_{j2}^h|-|R_{j2}^h|)\\ &\le
C(|P_{j1}^h|+|R_{j1}^h|+|P_{j2}^h|+|R_{j2}^h|)(|P_{j1}^h-R_{j1}^h|+|P_{j2}^h-R_{j2}^h|).
\end{align*}
Then, using a similar argument for $R_{3j}$ and H\"older's inequality it follows from \eqref{r_p_prime}  that
\begin{equation}
\label{r33}
\|R_{j3}^h\|_{L^2(S)}^2+\|R_{3j}^h\|_{L^2(S)}^2\le C \sqrt{\frac{E_2}{h^2}+E_1},\;\;j=1,2.
\end{equation}
Next we use that
\begin{align}
|\det R_h'-1| = & |\det R_h'-\det P_h'| \nonumber\\ = &|R_{11}^h R_{22}^h-R_{12}^h R_{21}^h-P_{11}^h P_{22}^h+P_{12}^h P_{21}^h| \nonumber \\ \le &|R_{11}^h R_{22}^h-P_{11}^h R_{22}^h|+|P_{11}^h R_{22}^h-P_{11}^h P_{22}^h|+|P_{12}^h R_{21}^h-R_{12}^h R_{21}^h|  +|P_{12}^h P_{21}^h-P_{12}^h R_{21}^h| \nonumber \\ \le & (3
|R_h|+|P_h'|)|R_h'-P_h'|. \label{det_r_p}
\end{align}
On the other hand, it holds
\begin{eqnarray*}
1=\det R_h= R_{33}^h \det R_h' -R_{32}^h(R_{11}^h R_{23}^h-R_{13}^hR_{21}^h) +R_{31}^h(R_{12}^h R_{23}^h-R_{13}^hR_{22}^h).
\end{eqnarray*}
Therefore,
\begin{align}
|R_{33}^h-1|^2\le  &C \left(|R_{33}^h|^2| \det R_h'-1|^2 +|R_{32}^h|^2(|R_{11}^h|^2+|R_{21}^h|^2) (|R_{23}^h|^2+|R_{13}^h|^2) \right.\nonumber \\ &\left. +|R_{31}^h|^2(|R_{12}^h|^2+|R_{22}^h|^2)(|R_{23}^h|^2+|R_{13}^h|^2) \right). \label{det_r}
\end{align}
Hence, integrating and using $|R_{33}^h|\leq 1$, \eqref{det_r_p}, \eqref{det_r} and \eqref{r33} we conclude
\begin{equation}
\label{det_r1}
\|R_{33}^h-1\|_{L^2(S)}^2\le C\left(\frac{E_2}{h^2}+E_1 \right),
\end{equation}
which concludes the proof of \eqref{r_p}.
\end{proof}
Let us now assume that $I^h(y^h)\le C E_h$, where
$$\lim\limits_{h\to 0}\frac{E_h}{h^2}=0.$$
Then it follows from the structure of the functional \eqref{func_i} and properties \ref{A3}, \ref{A4} that $E_2\le CE_h$ and $E_1\le C\frac{E_h}{h^2}$, where $E_i$ are defined in \eqref{e}.
In order to normalize the functions $y_h$, we prove the following lemma, c.f.~\cite{FGM_Der}.
\begin{lemma}
Let  $S\subseteq \mathbb R^2$ be a Lipschitz domain and $\Omega =S\times (-\frac{1}{2}, \frac{1}{2})$.  Let a sequence $y^h$ satisfy $I^h(y^h)\le C E_h$, where
\begin{equation}
\label{tend}
\lim\limits_{h\to 0}\frac{E_h}{h^2}=0.
\end{equation}
Let assumptions \ref{A3}, \ref{A4} be satisfied.
    Then there exist a rotation $B_h \in SO(3)$  and $c^h \in \R^3$ such that the following assertions hold for $B_h$ and the functions 
    \begin{equation}
\label{tilde_y}
\tilde y^h:=(B_h)^T(P_h)^T(y^h-c^h).
\end{equation}
The rotation $B_h$ has the form
    $B_h=\left(\begin{array}{ccc}B_h'& &0\\  & &0\\0&0&1    \end{array}\right)$, \nonumber\\
    \begin{eqnarray}
&&|B_h'-Id|\le Ch^{-1}\sqrt{E_h},\label{in_plane1}\\
&&\int\limits_\Omega (\tilde y_{1,2}^h-\tilde y_{2,1}^h) dx=0,\text{\qquad and}\label{in_plane}\\
&&\int\limits_\Omega \tilde y^h-\left(\begin{array}{c} x'\\ hx_3\end{array}\right) dx=0.\label{int_aver}
\end{eqnarray}
\label{lem:bh}
\end{lemma}
\begin{proof}
    We suppose that $h$ is small enough such that $\frac{E_h}{h}<1$. 
It is easy to see that \eqref{int_aver} is satisfied after a proper choice of the constants $c^h$. We use the notation
$$
\bar y^h=(P_h)^Ty^h.
$$
It follows from \eqref{q5} with H\"older's inequality and $E_2\leq CE_h$ that
\begin{equation}
\label{aver}
\left|\frac{1}{|\Omega|}\int\limits_\Omega\nabla'  {{\overline y'}^h} dx-Id\right|\le C \frac{\sqrt{E_h}}{h}.
\end{equation}
It remains to find $\theta \in (-\pi,\pi]$ such that the rotation ${B_h}'\in SO(2)$, 
$$
{B_h}'=\left(\begin{array}{cc}
   \cos\theta  &  -\sin\theta\\
    \sin\theta & \cos\theta
\end{array}  \right)
$$
satisfies \eqref{in_plane1} and \eqref{in_plane} .
It follows from \eqref{aver} that for $h>0$ small enough
\begin{eqnarray}
\label{below}
\left(\int\limits_\Omega ( {{\bar y}_{1,1}^h}+{{\bar y}_{2,2}^h})dx\right)^2\ge C\left(1-\left(\int\limits_\Omega ( {{\bar y}_{1,1}^h}-1)dx\right)^2-\left(\int\limits_\Omega ( {{\bar y}_{2,2}^h}-1)dx\right)^2\right)\ge C\left(1-\frac{E_h}{h^2}\right)>0.
\end{eqnarray}
In case 
$\int\limits_\Omega (\bar y_{1,2}^h- \bar y_{2,1}^h) dx=0$, 
we can choose ${B_h}'=Id$. If this is not the case, in order to satisfy property \eqref{in_plane}, we choose $\theta$ such that
$$
\cot \theta= \frac{\int\limits_\Omega ( {{\bar y}_{1,1}^h}+ {{\bar y}_{2,2}^h})dx}{\int\limits_\Omega ( {{\bar y}_{2,1}^h}- {{\bar y}_{1,2}^h})dx}.
$$
Furthermore, using \eqref{aver} and \eqref{below} we then obtain
\begin{eqnarray}
\label{sin}
|\sin\theta |=\sqrt{\frac{1}{1+\cot^2\theta}}\le C\frac{|\int\limits_\Omega ( {{\bar y}_{2,1}^h}- {{\bar y}_{1,2}^h})dx|}{\sqrt{\left(\int\limits_\Omega ( {{\bar y}_{2,1}^h}- {{\bar y}_{1,2}^h})dx\right)^2+\left( \int\limits_\Omega ( {{\bar y}_{1,1}^h}+ {{\bar y}_{2,2}^h})dx  \right)^2}}
\le C\frac{\sqrt{E_h}}{h}
\end{eqnarray}
and
\begin{eqnarray}
\label{cos}
|\cos\theta-1 |= \left|\frac{|\cot\theta|}{\sqrt{1+\cot^2\theta}}-1\right|\le C\frac{\left(\int\limits_\Omega ( {{\bar y}_{2,1}^h}- {{\bar y}_{1,2}^h})dx\right)^2}{{\left(\int\limits_\Omega ( {{\bar y}_{2,1}^h}- {{\bar y}_{1,2}^h})dx\right)^2+\left( \int\limits_\Omega ( {{\bar y}_{1,1}^h}+ {{\bar y}_{2,2}^h})dx  \right)^2}}\le C\frac{{E_h}}{h^2}.
\end{eqnarray}
Combining \eqref{sin} and \eqref{cos} we obtain \eqref{in_plane1}.

\end{proof}
Next, we define
\begin{eqnarray}
\tilde U^h(x')=\int\limits_I\left(\begin{array}{c}\tilde y_1^h\\\tilde y_2^h\end{array}\right)(x',x_3)-x' dx_3
\text{\qquad and\qquad}V^h(x')=\int\limits_I y_3^hdx_3.
\label{eq:tildeUV}
\end{eqnarray}
\begin{lemma}
\label{l_tilde_y}
Let  $S\subseteq \mathbb R^2$ be a Lipschitz domain and $\Omega =S\times (-\frac{1}{2}, \frac{1}{2})$.  Let $y^h$ be a sequence with $I^h(y^h)\le C E_h$, where $E_h$ satsifies \eqref{tend}. Let assumptions \ref{A3} and \ref{A4} be satisfied.
Then there exists $\tilde{u} \in W^{1,2}(S;\R^2)$ such that up to a (non-relabeled) subsequence
\begin{equation}
\label{tilde_u_conv}
\tilde u^h:=\min\left\{\frac{1}{\sqrt{E_h}}, \frac{h^2}{E_h}\right\}\tilde U_h\rightharpoonup \tilde u\quad \text{in} \;\; W^{1,2}(S;\mathbb R^2).
\end{equation}
\end{lemma}
\begin{proof}
We set $\tilde A_h':=\frac{h}{\sqrt{E_h}}((B_h')^T(P_h')^TQ_h-Id)$. Using  \ref{A4}, \eqref{q3}  \eqref{q6}, and \eqref{in_plane1} we obtain
$$\|\tilde A_h'\|_{W^{1,2}(S)}^2\le \frac{Ch^2}{E_h}(\|\nabla ' Q_h\|_{L^2(S)}^2+\|Q_h-P_h'\|_{L^2(S)}^2+\|B_h'-Id\|_{L^2(S)}^2)\le C,$$
and thus by the Sobolev embedding 
\begin{equation}
\label{a_prime_lp}
\|\tilde A_h'\|_{L^{p}(S)}\le C,\qquad \forall\;  1\le p<\infty.
\end{equation}
Since
\begin{eqnarray*}
&&(\tilde A_h')^T\tilde A_h'=\frac{h^2}{E_h}(Q_h^TP_h'B_h'-Id)((B_h')^T(P_h')^TQ_h-Id)=-\frac{h}{\sqrt{E_h}}((\tilde A_h')^T+\tilde A_h')\\&=&-\frac{h}{2\sqrt{E_h}}\sym  \tilde A_h'=-\sym \frac{h^2}{2E_h}((B_h')^T(P_h')^TQ_h-Id)
\end{eqnarray*}
we obtain from \eqref{a_prime_lp} with H\"older's inequality that
\begin{equation}
\label{a_prime}
\frac{h^2}{E_h}\|\sym ((B_h')^T(P_h')^TQ_h-Id)\|_{L^2(S)}\le C\|\tilde A_h'\|_{L^{4}(S)}^2\le C.
\end{equation}
Then, combining this with \eqref{q1}  we arrive at
\begin{align*}
&\|\sym  \nabla ' \tilde U_h\|_{L^{2}(S)} \\ =&\|\sym \left( \int\limits_I \nabla' (\tilde  y^h)' dx_3-Id \right)\|_{L^{2}(S)}\\ \le &\|\sym \left( \int\limits_I \nabla' (\tilde  y^h)' dx_3-(B_h')^T(P_h')^TQ_h \right)\|_{L^{2}(S)} +\|\sym ((B_h')^T(P_h')^TQ_h-Id)\|_{L^{2}(S)} \\
\le &C(\sqrt{E_h}+\frac{E_h}{h^2}).
\end{align*}
Using Korn's inequality and  normalizations \eqref{int_aver}, \eqref{in_plane} we infer the assertion of the lemma.
\end{proof}
Now we consider
\begin{align}
&U^h(x'):=\int\limits_I\left(\begin{array}{c} y_1^h\\y_2^h\end{array}\right)(x',x_3)-x' dx_3,\label{uu}
\end{align}
and recall that the quantities $P_h'$, $B_h'$ and $Q_h$ were introduced in Corollary \ref{cor1}, Lemma \ref{lem:bh}, and Theorem \ref{th_rig1}, respectively.
\begin{lemma}
\label{l_y}
Let a sequence $y^h$ satisfy the assumptions of Lemma \ref{l_tilde_y}.
Then it holds
\begin{align}
&|P_h'B_h'-Id|^2\le C \max\left(\frac{E_h^2}{h^4}, E_h\right),\label{p_conv}\\
&|P_h'-Id|^2\le C \frac{E_h}{h^2},\text{\qquad and}\label{p_conv1}\\
&\|Q_h(x')-Id\|_{L^2(S)}^2\le C \frac{E_h}{h^2}.\label{p_conv2}
\end{align}
Moreover, there exists $u \in W^{1,2}(S;\R^2)$ such that up to a (non-relabeled) subsequence 
\begin{equation}
\label{u_conv}
u^h:=\min\left(\frac{1}{\sqrt{E_h}}, \frac{h^2}{E_h}\right)U_h\rightharpoonup  u\quad \text{in} \;\; W^{1,2}(S;\mathbb R^2).
\end{equation}
\end{lemma}
\begin{proof}
The proof is similar to that presented in \cite{LM}. From \eqref{tilde_y} we infer the relation 
\begin{equation}
\label{yprime}
y_h'=P_h'B_h'(\tilde y^h)'+(c^h)'.
\end{equation}
It follows from \eqref{yprime} and the definitions of $u^h$  and $\tilde u^h$ that
\begin{eqnarray}
\label{u_tilde_u}
\max( h^{-2}E_h,\sqrt{E_h})u^h=(P_h'B_h'-Id)x'+\max( h^{-2}E_h,\sqrt{E_h})P_h'B_h'\tilde u^h+(c^h)'.
\end{eqnarray}
Due to the Dirichlet boundary conditions of functions with finite energy we have that $u^h=0$ on $\partial S$.
Combining \eqref{tilde_u_conv} and \eqref{u_tilde_u} and using the trace theorem we obtain
\begin{eqnarray}
\label{norm}
&&\|(P_h'B_h'-Id)x'+ (c^h)'\|_{L^2(\partial S)}^2\le \max( h^{-4}E_h^2, E_h)\|\tilde u^h\|_{L^2(\partial S)}^2\nonumber\\&\le& C\max( h^{-4}E_h^2, E_h)\|\tilde u^h\|_{W^{1,2}(S)}^2 \le  C\max( h^{-4}E_h^2, E_h).
\end{eqnarray}
Without loss of generality we may assume
$$\int\limits_{\partial S} x' \,d \mathcal{H}^1 =0\text{\qquad and\qquad}\int\limits_{\partial S} |x'|^2 \,d \mathcal{H}^1 >0.$$
For any $P\in SO(2)$ a straightforward computation shows, c.f.~\cite{LM},
\begin{equation}
\label{norm1}
2|(P-Id)x'|^2=|(P-Id)|^2|x'|^2.
\end{equation}
Therefore, \eqref{norm} yields 
\begin{equation}
\label{norm2}
|(c^h)'|\le C\max( h^{-2}E_h,\sqrt{E_h}).
\end{equation}
This together with \eqref{norm} and \eqref{norm1} gives \eqref{p_conv}.
Collecting \eqref{p_conv},\eqref{u_tilde_u}, \eqref{norm2} we arrive at
$$\|\nabla'u_h\|_{L^2( S)}^2\le C.$$
Taking into account the Dirichlet boundary conditions we obtain \eqref{u_conv}. 
Combining \eqref{in_plane1} and \eqref{p_conv} we get \eqref{p_conv1}. Finally, \eqref{p_conv2} is a consequence of \eqref{q6} and \eqref{p_conv1}.
\end{proof}

\begin{lemma}
\label{grad_bound}
Under the assumptions of Lemma \ref{l_y} and \ref{A3} it holds
\begin{align}
&\|\sym (R_h-Id)\|_{L^2(S)}^2\le C \frac{E_h}{h^2} \label{r_i}\text{\qquad and}\\
&\|\nabla_hy^h-Id\|_{L^2(\Omega)}^2\le C \frac{E_h}{h^2}\label{y_i}.
\end{align}
Moreover, there exists $v \in W^{1,2}(S)$ such that it holds up to a (non-relabeled) subsequence that 
\begin{equation}
\label{v_conv}
v^h(x'):=\frac{h}{\sqrt{E_h}} V^h(x') \rightharpoonup  v\quad \text{in} \;\; W^{1,2}(S).
\end{equation}
\end{lemma}
\begin{proof}
Since $R_h$ is a rotation we have for $i=1,2$
$$R_{i1}^h R_{31}^h+R_{i2}^hR_{32}^h+R_{i3}^hR_{33}^h=0.$$
Therefore,
\begin{align*}
&R_{21}^h R_{31}^h+(R_{22}^h-1)R_{32}^h+R_{23}^h(R_{33}^h-1)=R_{32}^h+R_{23}^h,\\
&(R_{11}^h-1) R_{31}^h+R_{12}^hR_{32}^h+R_{13}^h(R_{33}^h-1)=R_{31}^h+R_{13}^h.
\end{align*}
Consequently, using  \eqref{r_p_prime}, \eqref{r_p} and  \eqref{p_conv1}  we obtain
\begin{eqnarray}
\label{1}
\|R_{32}^h+R_{23}^h\|_{L^2(S)}^2  \le \|R_{21}^h\|_{L^2(S)}^2 \|R_{31}^h\|_{L^\infty(S)}^2+\|(R_{22}^h-1)\|_{L^2(S)}^2\|R_{32}^h\|_{L^\infty(S)}^2\nonumber\\+\|R_{23}^h\|_{L^\infty(S)}^2\|(R_{33}^h-1)\|_{L^2(S)}^2\le C\frac{E_h}{h^2},
\end{eqnarray}
\begin{eqnarray}
\label{2}
\|R_{31}^h+R_{13}^h\|_{L^2(S)}^2\le \|(R_{11}^h-1) \|_{L^2(S)}^2\|R_{31}^h\|_{L^\infty(S)}^2+\|R_{12}^h\|_{L^2(S)}^2\|R_{32}^h\|_{L^\infty(S)}^2\nonumber\\+\|R_{13}^h\|_{L^\infty(S)}^2\|(R_{33}^h-1)\|_{L^2(S)}^2\le C\frac{E_h}{h^2}. 
\end{eqnarray}
Then, it follows from \eqref{r_p_prime} \eqref{det_r1}, \eqref{p_conv1}, \eqref{1}, \eqref{2} that
\begin{eqnarray*}
\|\sym (R_h-Id)\|_{L^2(S)}^2\le C(\|R_h'-P_h'\|_{L^2(S)}^2+\|P_h'-I\|_{L^2(S)}^2\\+\|R_{31}^h+R_{13}^h\|_{L^2(S)}^2 +\|R_{32}^h+R_{23}^h\|_{L^2(S)}^2+\|(R_{33}^h-1)\|_{L^2(S)}^2) \le C\frac{E_h}{h^2}.
\end{eqnarray*}
Therefore, using \eqref{r1} we obtain
\begin{eqnarray*}
\|\sym (\nabla_h y^h-Id)\|_{L^2(\Omega)}^2\le C(\|\sym (R_h-Id)\|_{L^2(S)}^2\\+\|\nabla_h y^h-R_h\|_{L^2(\Omega)}^2)\le C\frac{E_h}{h^2}.
\end{eqnarray*}
Using the Dirichlet boundary conditions and Proposition 1 in \cite{FGM_Rig} we obtain \eqref{y_i}. Then \eqref{v_conv} is a consequence of \eqref{y_i}.
\end{proof}

\begin{theorem}
\label{th_rig2}
Suppose that $S\in \mathbb R^2$ is a Lipschitz domain and $\Omega =S\times (-\frac{1}{2}, \frac{1}{2})$.  Let $y^h$ be a sequence such that  $I^h(y^h)\le C E_h$, where $E_h$ satisfies $\lim\limits_{h\to 0}\frac{E_h}{h^2}=0.
$
Then  for any fixed $x_3\in I$ and $h$ small enough there exists 
a map
$ T_h(x', x_3):\Omega\to SO(2)$
such that  
\begin{align}
&\|\nabla' (y^h)'(x)-T_h(x)\|_{L^2(\Omega)}^2\le CE_h, 
\label{t1}\\
&\|T_h(x',x_3)-Id\|_{L^p(\Omega)}^2\le C\max\left(E_h^{\frac{2}{p}},\frac{E_h}{h^2}\right),\;\, 2\le p<\infty. \label{t2}
\end{align}
\end{theorem}
\begin{proof}
\begin{enumerate}
\item First we prove that there exists a map  $\tilde T_h(x', x_3):\Omega\to \mathbb R^{2\times 2}$ 
such that
\begin{align}
\|\nabla' (y^h)'(x)-\tilde T_h(x)\|_{L^2(\Omega)}^2\le CE_h, 
\label{t1_tilde}\\
\|\nabla' \tilde T_h(x',x_3)\|_{L^2(\Omega)}^2\le C\frac{E_h}{h^2}. \label{t2_tilde}
\end{align}

As in the proof of Theorem \ref{th_rig1} we consider an open subset $U \subseteq \mathbb R^2$ and $K\subseteq U$ compact such that $\dist_\infty (K, \partial U)>3h$. For each point $x'\in K$ we consider the square 
$$  S_{x',h}= x'+(0,h)^2$$
with lower left corner  $x'$ and
 define the map 
\begin{equation}
\label{def_t}
\hat T_{h}(x',x_3)= \int\limits_{S_{x',h}}h^{-2} \psi\left( \frac{x'-z'}{h}  \right)  F(z',x_3) dz',
\end{equation}
where $F$ is defined in \eqref{notatf} and $\psi$ is a standard mollifier.
We also use the notation $\psi_h(\cdot)=h^{-2}\psi(\cdot/h)$.
Using the Hölder's inequality and the rotation $Q_{x,x_3,h}$ defined in \eqref{s_rig} we get
\begin{align}
&\left|\hat T_h(x',x_3)- Q_{x',x_3, h} \right|^2\nonumber \\ = &\left|\int\limits_{S_{x',h}} h^{-2}\psi\left( \frac{x'-z'}{h}  \right) F(z', x_3)  dz'-\int\limits_{S_{x',h}} h^{-2}\psi\left( \frac{x'-z'}{h}  \right)Q_{x',x_3,h} dz' \right|^2\nonumber \\
\le &\frac{C}{h^2} \int\limits_{S_{x',h}} \left | F(z', x_3)-  Q_{x',x_3, h}\right|^2 dz'
\nonumber \\
\le &\frac{C}{h^2} \int\limits_{S_{x',h}} \dist^2(F(z',x_3), SO(2)) dz'. \label{hatr_rot_t}
\end{align}
Since $\int \nabla \psi_h=0$, using \eqref{s_i_rig} and H\"older's inequality we find for any point $\tilde x'\in S_{x',h}$
\begin{align}
&\left|\nabla' \hat T_h(\tilde x',x_3)\right|^2 \nonumber \\= &
\left|\int\limits_{S_{\tilde x',h}} h^{-3}(\nabla \psi)\left( \frac{\tilde x'-z'}{h}  \right) F(z', x_3)  dz'-\int\limits_{S_{\tilde x',h}} h^{-3}(\nabla \psi)\left( \frac{\tilde x'-z'}{h}  \right) Q_{x',x_3,h}  dz'  \right|^2 \nonumber \\
\le &\frac{C}{h^4} \int\limits_{S_{x',2h}} \left | F(z', x_3)-  Q_{x',x_3, h} \right|^2 dz'\nonumber \\
\le &\frac{C}{h^4} \int\limits_{S_{x',2h}} \dist^2(F(z',x_3), SO(2)) dz'. \label{hatr_nabla_t}
\end{align}
Integrating this inequality over $S_{x',h}\times I$ yields
\begin{equation}
\label{hatr_nabla2_t}
\int\limits_{S_{x',h}\times I}\left|\nabla' \hat T_h(z',x_3)\right|^2 dz'dx_3\le  \frac{C}{h^2} \int\limits_{S_{x',h}\times I} \dist^2(F(z',x_3), SO(2)) dz' dx_3 .
\end{equation}
For any point $\tilde x'\in S_{x',h}$ we have
\begin{align}
&\left|\hat T_h(x',x_3)- \hat T_h(\tilde x',x_3) \right|^2 \nonumber \\\le &\left|\int\limits_{S_{x',h}} h^{-2}\psi\left( \frac{x'-z'}{h}  \right) F(z', x_3) dz'-\int\limits_{S_{\tilde x',h}} h^{-2}\psi\left( \frac{\tilde x'-z'}{h}  \right) F(z', x_3)  dz'  \right|^2 \nonumber \\ \le 
&\left|\int\limits_{S_{x',h}} h^{-2}\psi\left( \frac{x'-z'}{h}  \right) F(z', x_3)  dz'-\int\limits_{S_{ x',h}} h^{-2}\psi\left( \frac{ x'-z'}{h}  \right) Q_{x',x_3, h}  dz'  \right|^2 \nonumber \\&+\left|\int\limits_{S_{\tilde x',h}} h^{-2}\psi\left( \frac{\tilde x'-z'}{h}  \right) F(z', x_3)  dz'-\int\limits_{S_{ \tilde x',h}} h^{-2}\psi\left( \frac{ \tilde x'-z'}{h}  \right) Q_{x',x_3, h}  dz'  \right|^2 \nonumber \\
\le &\frac{C}{h^2} \int\limits_{S_{x',2h}} \left | F(z', x_3)-  Q_{x',x_3, h} \right|^2 dz' \nonumber \\
\le &\frac{C}{h^2} \int\limits_{S_{x',2h}} \dist^2(F(z',x_3), SO(2)) dz'. \label{hatr_hatr_t}
\end{align}
Combining this with \eqref{s_rig} and \eqref{hatr_rot_t} yields
\begin{align}
&\int\limits_{S_{x',h}}\left|\hat T_h(z',x_3)- F(z', z_3) \right|^2dz' \nonumber \\
\le &\int\limits_{S_{x',h}}\left|\hat T_h(z',x_3)- \hat T_h(x',x_3)\right|^2dz' +
\int\limits_{S_{x',h}}\left|\hat T_h(x',x_3)-Q_{x',x_3,h}  \right|^2dz'+\int\limits_{S_{x',h}}\left| Q_{x',x_3,h}  - F(z', x_3) \right|^2dz' \nonumber\\
\le &C \int\limits_{S_{x',h}} \dist^2(F(z',x_3), SO(2)) dz' \label{hatr_grad_t}
\end{align}
Finally, setting $g(\zeta)=\dist (\hat T_h(x'+h\zeta,x_3), SO(2))$ and taking into account \eqref{hatr_nabla_t},  \eqref{hatr_grad_t}
we come to 
\begin{align*}
&\int_{(0,1)^2}|g|^2 d\zeta+\sup\limits_{(0,1)^2}|\nabla' g|^2 \, dz'\\
=&\int\limits_{S_{x',h}} \dist^2(\hat T_h(z',x_3), SO(2)) dz'+\sup\limits_{S_{x',h}}|\nabla' \dist(\hat T_h(z',x_3), SO(2)) |^2\\
\le &C\bigg( \int\limits_{S_{x',h}} |\hat T_h(z',x_3)-  
 F(z',x_3)|^2dz'+\int\limits_{S_{x',h}} \dist^2(
 F(z',x_3), SO(2)) dz' +\sup\limits_{S_{x',h}}|\nabla' \hat T_h(z',x_3)|^2 \bigg) \\
\le&\frac{C}{h^2} \int\limits_{S_{x',2h}} \dist^2(F(z',x_3), SO(2)) dz'.
\end{align*}
Consequently,
\begin{eqnarray*}
\dist^2 (\hat T_h(\tilde x',x_3), SO(2))&\le& \sup\limits_{(0,1)^2}|g|^2\le C(\int_{(0,1)^2}|g|^2 d\zeta+\sup\limits_{(0,1)^2}|\nabla' g|^2)\\
&\le&\frac{C}{h^2} \int\limits_{S_{x',2h}} \dist^2(F(z',x_3), SO(2)) dz',
\end{eqnarray*}
and thus
\begin{eqnarray}
\label{hatr_bound2_t}
\int\limits_{S_{x',h}}\dist^2 (\hat T_h(\tilde x',x_3), SO(2)) d \tilde x'
&\le&\frac{C}{h^2} \int\limits_{S_{x',h}}\int\limits_{S_{\tilde x',2h}} \dist^2(F(z',x_3), SO(2)) dz' d\tilde x'\nonumber\\
&\le&\frac{C}{h^2} \int\limits_{S_{x',h}}\int\limits_{S_{x',4h}} \dist^2(F(z',x_3), SO(2)) dz' d\tilde x'\nonumber\\&\le&
C \int\limits_{S_{x',h}} \dist^2(F(z',x_3), SO(2)) dz'.
\end{eqnarray}
Now we consider a square lattice of size $h$ in  $\mathbb R^2$ and sum the inequalities \eqref{hatr_nabla2_t}, \eqref{hatr_grad_t} over all squares which intersect $K$. This yields
\begin{align}
&\int\limits_{K\times I}\left(|\hat T_h(x',x_3)- 
 F(x',x_3)|^2+h^2|\nabla' \hat T_h(x',x_3)|^2\right)dx'dx_3 \le C \int\limits_{U\times I} \dist^2(F(x), SO(2)) dx.
\label{hatr_k_t}
\end{align}
Following \cite{FGM_Der} we consider $S$ to be locally the epigraph of a Lipschitz function, then arguing as in \cite{FGM_Der} and applying the above estimates we get
\begin{eqnarray}
\label{hatr_k_st}
\int\limits_{K\cap S\times I}\left(|\hat T_h(x',x_3)-
 F(x',x_3)|^2+h^2|\nabla' \hat T_h(x',x_3)|^2\right)dx'dx_3\le C \int\limits_{U\times I} \dist^2(F(x), SO(2)) dx.
\end{eqnarray}
Since $S$ is Lipschitz $\bar S$ can be covered by a finite number of open sets $U_0,..., U_l$, where $\bar U_0\subset S$. Denote $\hat T_i$ the maps constructed in previous steps for $T_i$ and consider a partition of unity correspondent to the cover $\{T_i\}$
$$ \eta_i\in C_0^\infty(U_i), \;\qquad\eta_i\ge 0, \;\qquad\sum\limits_{i=0}^l\eta_i(x')=1,\,\forall x'\in S.$$
If we set
\begin{equation}
\label{tilde_rt}
\tilde T_h=\sum\limits_{i=0}^l\eta_i \hat T_i
\end{equation}
then 
\begin{align*}
&\tilde T_h- F =\sum\limits_{i=0}^l\eta_i (\hat T_i- F ),\text{\qquad and\qquad}
\nabla'\tilde T_h=\sum\limits_{i=0}^l\eta_i \nabla'\hat T_i+\sum\limits_{i=0}^l\nabla'\eta_i (\hat T_i-F ).
\end{align*}
Consequently, applying \eqref{hatr_k_t} and \eqref{hatr_k_st} with $K_i=\text{supp} \,\eta_i$ we obtain
\begin{equation}
\label{tildert}
\int\limits_{\Omega}\left(|\tilde T_h(x',x_3)-
 F(x',x_3)|^2+h^2|\nabla' \tilde T_h(x',x_3)|^2\right)dx\le C E_h,
\end{equation}
which proves \eqref{t1_tilde} and \eqref{t2_tilde}.
It follows from \eqref{hatr_bound2_t} that
\begin{equation*}
\int\limits_{\Omega}\dist^2 (\tilde T_h( x',x_3), SO(2)) d  x
\le CE_h,
\end{equation*}
and hence
\begin{eqnarray}
\label{tilder_bound2_t}
\int\limits_{\Omega}|\tilde T_h( x',x_3)|^2 dx\le C\int\limits_{\Omega} (\inf\limits_{P\in SO(2)}|\tilde T_h( x',x_3)-P|^2+ 1)\, dx \le C(E_h+1)\le C.
\end{eqnarray}
\item Now the rotation $T_h$ can be obtained by projecting $\tilde T_h$ onto $SO(2)$. Since $SO(2)$ is a smooth manifold, there exists a tubular neighborhood $\mathcal U$ of $SO(2)$ 
such that the projection $\pi:\mathcal U\to SO(2)$ is smooth. Note that there exists a $\delta>0$ such that for all $M\not\in \mathcal{U}$ there holds $\dist(M,SO(2))\geq \delta$. We define
$$
T_h(x',x_3)=\tilde \pi(\tilde T_h(x',x_3))=\begin{cases}{c}\pi(\tilde T_h(x',x_3)) &\text{ if } \tilde T_h(x',x_3)\in \mathcal U,\\  Id &\text{ if } \tilde T_h(x',x_3)\notin \mathcal U.\end{cases}
$$
Then in case $\tilde T_h(x',x_3)\in\mathcal U$ we have 
$|T_h(x',x_3)-\tilde T_h(x',x_3)|= \dist(\tilde T_h(x',x_3), SO(2))$.
If $\tilde T_h(x',x_3)\notin\mathcal U$ then
$\dist^2 (\tilde T_h(x',x_3), SO(2))dx>\delta^2$, and hence
\begin{align*}
|T_h(x',x_3)-\tilde T_h(x',x_3)|^2 = &|Id-\tilde T_h(x',x_3)|^2 \\
\leq &C( 1 + \operatorname{dist}^2(T_h(x',x_3),SO(2)) ) \\ \leq &C \left( \frac{1}{\delta^2} + 1 \right) \operatorname{dist}^2(T_h(x',x_3),SO(2)).
\end{align*}
Integration over $\Omega$ then yields
$$
\int\limits_{\Omega}|T_h(x',x_3)-\tilde T_h(x',x_3)|^2 dx\le C\int\limits_\Omega \dist^2(\tilde T_h(x',x_3), SO(2)) dx\le C E_h,
$$
which together with \eqref{t1_tilde}  gives \eqref{t1}.
\item 
 Let $\Omega'$ be any compact subset of $\Omega$ and $|s|<\dist(\Omega',\partial\Omega)$. It follows from \eqref{def_t}, integration by parts and Hölder's inequality that for any $s>0$
\begin{align}
&\left| \frac{\hat T_h(x', x_3+s)-\hat T_h(x', x_3)}{s} \right|^2 \nonumber\\=&\left|\int\limits_{S_{x',h}} h^{-2}\psi\left( \frac{x'-z'}{h}\right)\frac{\nabla' (y^h)'(z', x_3+s)-\nabla' (y^h)'(z', x_3)}{s} dz'\right|^2 \nonumber \\
=&\left|\int\limits_{S_{x',h}} h^{-3}\nabla'\psi\left(\frac{x'-z'}{h}\right)\frac{ (y^h)'(z', x_3+s)- (y^h)'(z', x_3)}{s} dz'\right|^2 \nonumber \\
\le &Ch^{-4}\int\limits_{S_{x',h}} \left|\frac{(y^h)'(z', x_3+s)-(y^h)'(z', x_3)}{s}\right|^2 dz'. \label{t_diff}
\end{align}
After integration of \eqref{t_diff} over $S_{x',h}\times I$ we obtain
\begin{align*}
&\int\limits_{S_{x',h}\times I}\left| \frac{\hat T_h(x', x_3+s)-\hat T_h(x', x_3)}{s} \right|^2 dx\le Ch^{-2}\int\limits_{S_{x',h}\times I} \left|\frac{(y^h)'(x', x_3+s)-(y^h)'(x', x_3)}{s}\right|^2 dx.
\end{align*}
Therefore, since $\Omega'$ is arbitrary, we get for the map $\tilde T_h$ defined in \eqref{tilde_rt} that 
\begin{eqnarray*}
\left\| \frac{\tilde T_h(x', x_3+s)-\tilde T_h(x', x_3)}{s}\right\|_{L^2(\Omega)}^2\le \frac{C}{h^2}\left\|\frac{(y^h)'(x', x_3+s)-(y^h)'(x', x_3)}{s}\right\|_{L^2(\Omega)}^2.
\end{eqnarray*}
Passing to the limit $s\to 0$ and using \eqref{y_i} we obtain
\begin{equation*}
\|\partial_3 \tilde T_{h}(x',x_3)\|_{L^2(\Omega)}^2\le \frac{C}{h^2}\|\partial_3 (y^h)'(x', x_3)\|_{L^2(\Omega)}^2\le C\frac{E_h}{h^2}. 
\end{equation*}
Combining this with \eqref{tildert} yields
\begin{equation}
\label{nabla_t}
\|\nabla \tilde T_{h}(x',x_3)\|_{L^2(\Omega)}^2\le C\frac{E_h}{h^2}.
\end{equation}
This implies that for the constant map $S_h'=\frac{1}{|\Omega|}\int\limits_\Omega \tilde T_{h}(x) dx$ we have
\begin{align}
\|\int\limits_I \tilde T_{h}(x',x_3) d x_3-S_h'\|_{L^2(S)}^2  \le &C \|\tilde T_{h}(x',x_3)-S_h'\|_{L^2(\Omega)}^2 \le C\|\nabla \tilde T_{h}(x',x_3)\|_{L^2(\Omega)}^2\le  C\frac{E_h}{h^2}.
\label{nabla_ts}
\end{align}
Together with \eqref{q1}, \eqref{q6} and \eqref{t1_tilde} one may then derive 
\begin{align}
&|S_h'-P_h'|^2 \nonumber \\ \le &C(\|S_h'-\int\limits_I \tilde T_{h}(x',x_3) d x_3\|_{L^2(S)}^2+\|\int\limits_I \tilde T_{h}(x',x_3) d x_3-\int\limits_I\nabla' (y^h)'(x) dx_3\|_{L^2(S)}^2 \nonumber \\ &+ \|\int\limits_I \nabla' (y^h)'(x)  d x_3-Q_h(x')\|_{L^2(S)}^2+\|Q_h(x')-P_h'\|_{L^2(S)}^2) \nonumber \\ 
\le &C \frac{E_h}{h^2}. \label{sp}
\end{align}
Combining this with \eqref{p_conv1}, \eqref{nabla_ts} implies
\begin{align}
\label{tid0}
&\|\tilde T_{h}(x',x_3)-Id\|_{L^2(\Omega)}^2  \le C(\|\tilde T_{h}(x',x_3)-S_h'\|_{L^2(\Omega)}^2+|S_h'-P_h'|^2+|P_h'-Id|^2)\le C\frac{E_h}{h^2}.
\end{align}
By the Sobolev embedding theorem one may then derive from \eqref{nabla_t} and \eqref{tid0} for any $1\le p<\infty $ that
\begin{equation}
\label{tid1}
\|\tilde T_{h}(x',x_3)-Id\|_{L^p(\Omega)}^2\\\le C\|\tilde T_{h}(x',x_3)-Id\|_{W^{1,2}(\Omega)}^2\le C\frac{E_h}{h^2}.
\end{equation}
In case $\tilde T_{h}(x',x_3)\notin U$, it follows from the definition of $T_{h}(x',x_3)$ that  $\|T_{h}(x',x_3)-Id\|_{L^p(S)}^2=0$. In case $\tilde T_{h}(x',x_3)\in U$ it holds  $$|T_{h}(x',x_3) - \tilde{T}_{h}(x',x_3)| = \dist(\tilde T_h(x',x_3), SO(2))<\delta.$$
Hence, it holds for all $2 \leq p < \infty$ that 
\[
|T_{h}(x',x_3)-Id|^p \leq C \left( \dist^p(\tilde T_h(x',x_3), SO(2)) + |\tilde T(x',x_3)-Id|^p  \right).
\]
Consequently, it follows
for any $2\le p<\infty$ and any $x_3\in I$
\begin{equation*}
\|T_{h}(x',x_3)-Id\|_{L^p(S)}^p \leq C \left( \int\limits_S\dist^p(\tilde T_h(x',x_3), SO(2))dx' + \|\tilde T(x',x_3)-Id\|_{L^p(S)}^p  \right).
\end{equation*}
Integrating this inequality over $I$ with respect to $x_3$  we infer
\begin{align}
&\|T_{h}(x',x_3)-Id\|_{L^p(\Omega)}^2 \nonumber \\ \le & C \left(\int\limits_\Omega \dist^p(\tilde T_h(x',x_3), SO(2)) dx\right)^{\frac{2}{p}}+ C\|\tilde T_{h}(x',x_3)-Id\|_{L^p(\Omega)}^2 \nonumber\\\le &C\max\left(E_h^{\frac{2}{p}},\frac{E_h}{h^2}\right). \label{tid2}
\end{align}
\end{enumerate}
\end{proof}

\section{Dimension reduction}\label{sec:dimred}
We introduce the rotation 
\begin{equation}
\label{L_h}
L_h=\left(\begin{array}{ccc}Q_h& &0\\  & &0\\0&0&1    \end{array}\right).
\end{equation}
\begin{lemma}
\label{lem_g1}
Let $S\in \mathbb R^2$ be a Lipschitz domain and $\Omega =S\times (-\frac{1}{2}, \frac{1}{2})$.  Let $y^h$ be a sequence such that $I^h(y^h)\le C E_h$, where
\begin{equation}
\label{tend1}
\lim\limits_{h\to 0}\frac{E_h}{h^4}=1
\end{equation}
or
\begin{equation}
\label{tend2}
\lim\limits_{h\to 0}\frac{E_h}{h^4}=0.
\end{equation}
Then there exists $G \in L^2(\Omega;\R^{3\times 3})$ such that it holds up to a non-relabeled subsequence  
\begin{equation}
\label{gl2}
G_h:=\frac{L_h^T\nabla_h y^h -Id}{(E_h)^{1/2}/h}\rightharpoonup G \quad \text{in}\;\;L^2(\Omega;\mathbb R^{3\times 3}) \text{ as $h \to 0$}.
\end{equation}
Moreover, $G$ has the form
\[
G = \begin{pmatrix} 0 && 0 && G_{31} \\ 0 && 0 && G_{32} \\ \partial_1v && \partial_2 v && G_{33} \end{pmatrix}
\]
where  $G_{31}$ and $G_{32}$ do not depend on $x_3$.
\end{lemma}
\begin{proof}
1. For $h$ small enough and under assumption  \eqref{tend1} or \eqref{tend2} one can infer from \eqref{r1}, \eqref{q6}, \eqref{r_p_prime}, and \eqref{y_i} that
\begin{eqnarray*}
\|G_h\|_{L^2(\Omega)}^2\le C \left(\left\|\frac{\nabla' (y^h)'-Q_h}{(E_h)^{1/2}/h}\right\|_{L^2(\Omega)}^2+\left\| \frac{\nabla ' y_3}{(E_h)^{1/2}/h}\right\|_{L^2(\Omega)}^2+\left\| \frac{\partial_3 y'}{(E_h)^{1/2}}\right\|_{L^2(\Omega)}^2+\left\| \frac{\frac{\partial_3 y_3}{h}-1}{(E_h)^{1/2}/h}\right\|_{L^2(\Omega)}^2\right)\le C.
\end{eqnarray*}
   Therefore, \eqref{gl2} holds. \\[.2cm]
2. Let $\Omega'$ be any compact subset of $\Omega$ and $|s|<\dist(\Omega',\partial\Omega)$ and consider the difference quotients 
 $$H_h(x',x_3):=s^{-1}(G_h(x',x_3+s)-G_h(x',x_3)).$$
 Due to \eqref{gl2} 
 \begin{equation}
\label{hl2}
H_h\rightharpoonup H=s^{-1}(G(x',x_3+s)-G(x',x_3)) \quad \text{in}\;\;L^2(\Omega';\mathbb R^{3\times 3}).
\end{equation}
On the other hand,
$$H_h=L_h^T\frac{\nabla_h y^h(x',x_3+s)-\nabla_h y^h(x',x_3)}{s(E_h)^{1/2}/h}.$$
It follows from \eqref{q3} and \eqref{p_conv2} combined with the Sobolev embedding theorem that up to taking a subsequence it holds for all $1\leq p < \infty$ that  
\begin{equation}
    \label{lh_conv}
 L_h   \to  Id \quad \text{in}\;\;L^p(S;\mathbb R^{3\times 3}).
\end{equation}
Therefore,
 \begin{equation}
\label{lhl2}
L_hH_h=\frac{\nabla_h y^h(x',x_3+s)-\nabla_h y^h(x',x_3)}{s(E_h)^{1/2}/h}\rightharpoonup H \: \text{in}\;\;L^1(\Omega';\mathbb R^{3\times 3}).
\end{equation}
On the other hand,
\begin{eqnarray}
\label{lh}
L_hH_h&=&\frac{\nabla_h y^h(x',x_3+s)-\nabla_h y^h(x',x_3)}{s(E_h)^{1/2}/h}\nonumber\\
&=&\left(\frac{h^2}{(E_h)^{1/2}}\nabla'\frac{1}{s}\int\limits_{x_3}^{x_3+s}\frac{1}{h} \partial_3 y^h(x',z)  dz  \;\bigg| \;\; \frac{1}{(E_h)^{1/2}}\frac{1}{s}\int\limits_{x_3}^{x_3+s} \partial_{3}^2 y^h(x',z)  dz \right).
\end{eqnarray}
By definition of $I^h$ it holds $\int_{\Omega} \frac1{h^{\eps}} |\partial_3^2 (y^h)'|^2 \, dz \leq C E_h$ and thus  
$$\frac{1}{(E_h)^{1/2}}\partial_{3}^2 y'^h(x',z) \to (0,0)^T\quad \text{in}\;L^2(\Omega';\mathbb R^2). $$
Next, using \eqref{y_i} one deduces that
$$ \frac{h^2}{(E_h)^{1/2}}\nabla'\frac{1}{h} \partial_3 y^h(x',z)   \to \left( \begin{array}{cc}
     0&0  \\
     0&0 \\
     0&0
\end{array}\right)  \quad \text{in}\;W^{-1,2}(\Omega';\mathbb R^{2\times 3}). $$
Therefore, since $\Omega'$ is arbitrary it follows from \eqref{lh_conv},\eqref{lhl2} and \eqref{lh} that all the entries of $H$ are equal to zero, except maybe $H_{33}$, consequently, the  entries of  $G(x',x_3)$ do not depend on $x_3$, except maybe $G_{33}$. \\[.2cm]
3. Now we notice that $G_h'=\frac{Q_h^T\nabla' (y^h)' -Id}{(E_h)^{1/2}/h}$ and that \eqref{gl2} yields
\begin{equation}
\label{aver_g}
\int\limits_I G'_h(x',x_3) dx_3\rightharpoonup \int\limits_I G '(x')dx_3=G'(x'), \quad \text{in}\;\;L^2(S;\mathbb R^{2\times 2}).
\end{equation}
Then, due to \eqref{q1}
\begin{eqnarray*}\int\limits_I G_h'(x',x_3) dx_3=\int\limits_I \frac{Q_h^T\nabla' (y^h)' -Id}{(E_h)^{1/2}/h} dx_3=\frac{Q_h^T\int\limits_I\nabla' (y^h)'dx_3 -Id}{(E_h)^{1/2}/h}\to 0  \quad \text{in}\;\;L^2(S;\mathbb R^{2\times 2}),
\end{eqnarray*}
and therefore, 
\[
G' = \begin{pmatrix} 0 && 0 \\ 0 && 0 \end{pmatrix}.  
\]
Next, it holds for $\alpha=1,2$
$$\int\limits_IG_{3\alpha}^h dx_3=h\int\limits_I\partial_\alpha y_3/(E_h)^{1/2}dx_3=\partial_\alpha v_h\rightharpoonup \partial_\alpha v\quad \text{in}\;\;L^2(S),$$
which together with \eqref{aver_g} gives $G_{3\alpha}=\partial_\alpha v(x')$.
\end{proof}
\begin{remark}
\label{phi_def}
Let us denote
\begin{equation*}
\varphi_\alpha^h(x',x_3)=\frac{1}{(E_h)^{1/2}}\left(Q_h^T\partial_3 y'(x',x_3)\right)_\alpha\rightharpoonup \varphi_\alpha(x')=G_{\alpha 3} \quad \text{in}\;\;L^2(S;\mathbb R^{3\times 3}).
\end{equation*}
Then
$$G=\left( \begin{array}{ccc}
     0&0&\varphi_1(x')  \\
     0&0&\varphi_2(x')  \\
     \partial_1 v(x')&\partial_2 v(x')&G_{33}(x)
\end{array}\right).$$
\end{remark}

\begin{lemma}
\label{lem_g2} Let the assumptions of Lemma \ref{lem_g1} be satisfied.
Then  there exists $F\in L^2(\Omega;\R^{3\times 3}{)}$ such that it holds up to a non-relabeled subsequence
\begin{equation}
\label{fl2}
F_h:=\frac{Q_h^T\int\limits_I\nabla' (y^h)' dx_3-Id}{(E_h)^{1/2}}\rightharpoonup F(x') \quad \text{in}\;\;L^2(\Omega;\mathbb R^{3\times 3})
\end{equation}
as $h\to 0$
and
\begin{equation}
\label{fl22}
\sym  F_h\rightharpoonup \sym \nabla'u(x')\quad \text{in}\;\;L^2(\Omega;\mathbb R^{3\times 3}),
\end{equation}
 where $u$ is defined in \eqref{u_conv}.
\end{lemma}
\begin{proof}
It follows from \eqref{q1} that
   $$
\|F_h\|_{L^2(\Omega)}^2\le C \left\|\frac{\int\limits_I\nabla' (y^h)' dx_3-Q_h}{(E_h)^{1/2}}\right\|_{L^2(\Omega)}^2\le C,
$$
which implies\eqref{fl2}.\\[.2cm]
Now we rewrite $F_h$ as
   \begin{eqnarray}
 \label{f_dec}
F_h=\frac{\int\limits_I\nabla' (y^h)' dx_3-Id}{(E_h)^{1/2}} -\frac{Q_h-Id}{(E_h)^{1/2}}+(Q_h-Id)^T\frac{\int\limits_I\nabla' (y^h)' dx_3-Q_h}{(E_h)^{1/2}}
 \end{eqnarray}
 and set $A_h:=\frac{Q_h-Id}{(E_h)^{1/4}}$.
Using that $Q_h \in SO(2)$ we observe that
\begin{eqnarray}
  \label{sym_q}
  A_h^T A_h=\frac{(Q_h-Id)^T}{(E_h)^{1/4}}\frac{(Q_h-Id)}{(E_h)^{1/4}}=-\frac{2}{(E_h)^{1/4}}\sym A_h=-2\sym \frac{(Q_h-Id)}{(E_h)^{1/2}}.
   \end{eqnarray}
   It follows from \eqref{q3}, \eqref{p_conv2} and the Sobolev embedding that
  \begin{equation}
  \label{sym_q_bound}
A_h\to A \quad \text{in}\;\; L^p(S;\mathbb R^{2\times 2}), \quad 1\le p<\infty.
\end{equation}
   Combining \eqref{sym_q_bound} with \eqref{sym_q} we get
   \begin{equation}
   \label{sym_q_bound1}
\sym A_h\to 0\quad \text{in}\;\; L^p(S;\mathbb R^{2\times 2}) \quad 1\le p<\infty.  
\end{equation}
Therefore, it follows from \eqref{sym_q_bound} and \eqref{sym_q_bound1} that $\sym A=0$. It remains to identify $A_{12}$.
It follows from \eqref{q1}, \eqref{p_conv2} and \eqref{u_conv} that
\begin{equation}
  \label{sym_q_bound3}
A_{12}=\lim\limits_{h\to 0} (E_h)^{1/4}u_h=0\quad\text{in}\;\; L^p(S) \quad 1\le p<\infty, 
\end{equation}
which implies that $A=0$.
Then, \eqref{sym_q_bound} yields 
\begin{equation}
  \label{sym_q_bound4}
\sym \frac{(Q_h-Id)}{(E_h)^{1/2}}\to -\frac{A^2}{2}=0\quad\text{in}\; L^p(S;\mathbb R^{2\times 2}) \quad 1\le p<\infty. 
\end{equation}
Now we estimate the third term in \eqref{f_dec}. By \eqref{q1} and \eqref{p_conv2}  we have
\begin{eqnarray*}
\left\|(Q_h-Id)^T\frac{\int\limits_I\nabla' (y^h)' dx_3-Q_h}{(E_h)^{1/2}}\right\|_{L_1(\Omega)}\le \left\|Q_h-Id\right\|_{L^2(\Omega)} \left\|\frac{\int\limits_I\nabla' (y^h)' dx_3-Q_h}{(E_h)^{1/2}}\right\|_{L^2(\Omega)} \le C\frac{\sqrt{E_h}}{h}\to 0,
 \end{eqnarray*}
which together with \eqref{u_conv}, \eqref{sym_q_bound4} and \eqref{fl2} yields \eqref{fl22}. 
\end{proof}
\begin{lemma}
\label{lem_g3} Let the assumptions of Lemma \ref{lem_g1} be satisfied. Define the function $K \in L^2(\Omega;\R^{2\times 2})$ as $K(x',x_3) := x_3 \nabla' \varphi(x')$ where $\varphi$ is defined in \eqref{phi_def}.
Then  it holds up to to a non-relabeled subsequence as $h \to 0$
\begin{equation}
\label{kl2}
K_h:=\frac{Q_h^T(\nabla' (y^h)'-\int\limits_I\nabla' (y^h)' dx_3)}{(E_h)^{1/2}}\rightharpoonup K(x',x_3) \quad \text{in}\;\;L^2(\Omega;\mathbb R^{2\times 2}).
\end{equation}
\end{lemma}
\begin{proof}
We use the Dirichlet boundary condition in combination with Korn's inequality applied in the two-dimensional domains $S \times \{x_3\}$ to estimate
   \begin{align}
\|K_h\|\le &C\left\|\frac{\nabla' (y^h)'-\int\limits_I\nabla' (y^h)' dx_3}{(E_h)^{1/2}}\right\|_{L^{2}(\Omega)}^2 \nonumber \\ \le &C\left\|\frac{\sym (\nabla' (y^h)'-\int\limits_I\nabla' (y^h)' dx_3)}{(E_h)^{1/2}}\right\|_{L^{2}(\Omega)}^2\nonumber \\ \le
&C\left\|\frac{\nabla' (y^h)'-T_h}{(E_h)^{1/2}}\right\|_{L^{2}(\Omega)}^2+C\left\|\frac{\sym (T_h-\int\limits_IT_h dx_3)}{(E_h)^{1/2}}\right\|_{L^{2}(\Omega)}^2 \nonumber \\\le
&C\left\|\frac{\nabla' (y^h)'-T_h}{(E_h)^{1/2}}\right\|_{L^{2}(\Omega)}^2+C\left\|\frac{\sym (T_h-Id)}{(E_h)^{1/2}}\right\|_{L^{2}(\Omega)}^2, \label{tt}
 \end{align}
 where $T_h$ is the function from Theorem \ref{th_rig2}. The first term on the right-hand side is bounded due to \eqref{t1}. Therfore, it is left to show that also the second term is bounded. 
    We set
    $C_h':=\frac{T_h(x',x_3)-Id}{(E_h)^{1/4}}$ and compute
\begin{eqnarray*}
     (C_h')^T C_h'&=&\frac{T_h(x',x_3)^T-Id}{(E_h)^{1/4}}\cdot\frac{T_h(x',x_3)-Id}{(E_h)^{1/4}}\\&=&-\frac{1}{(E_h)^{1/4}}\left(\frac{T_h(x',x_3)^T-Id}{(E_h)^{1/4}}+\frac{T_h(x',x_3)-Id}{(E_h)^{1/4}}\right)\\&=&-2\frac{1}{(E_h)^{1/4}}\sym C_h'.
 \end{eqnarray*}
 Hence, we may estimate using \eqref{t2} and H\"older's inequality
 \[
 \left\|\frac{\sym (T_h-Id)}{(E_h)^{1/2}}\right\|_{L^{2}(\Omega)}^2 \leq \left\|\frac{\sym (T_h-Id)}{(E_h)^{1/4}}\right\|_{L^{4}(\Omega)}^4 \leq C \max\left\{ 1, \frac{E_h}{h^4} \right\} \leq C.
 \]
 \item 
  Let now again $\Omega'$ be any compact subset of $\Omega$ and $|s|<\dist(\Omega',\partial\Omega)$.
 We introduce the difference quotients 
 $$M_h(x',x_3):=s^{-1}(K_h(x',x_3+s)-K_h(x',x_3)).$$
 Due to \eqref{kl2} 
 \begin{equation}
\label{ml2}
M_h\rightharpoonup M=s^{-1}(K(x',x_3+s)-K(x',x_3)) \quad \text{in}\;\;L^2(\Omega';\mathbb R^{2\times 2}).
\end{equation}
On the other hand, we have by Remark \ref{phi_def}
\begin{align*}M_h & =Q_h^T\frac{\nabla' (y^h)'(x',x_3+s)-\nabla' (y^h)'(x',x_3)}{s(E_h)^{1/2}}\\ &=\frac{1}{(E_h)^{1/2}}\nabla'\frac{1}{s}Q_h^T\int\limits_{x_3}^{x_3+s}\partial_3 (y^h)'(x',z)  dz \rightharpoonup \nabla ' \varphi(x') \quad \text{in}\;W^{-1,2}(\Omega';\mathbb R^{2\times 2}).
\end{align*}
Since $\Omega'$ is arbitrary, the claimed form of $K$ follows.
\end{proof}
Now we prove our main result.\\
\begin{proof}[\bf Proof of Theorem \ref{th_conv1}]
\begin{enumerate}
\item[(i)] 
Let  $E_h=h^\sigma$. It follows from Theorem \ref{th_rig1}, Lemmas \ref{l_y}, \ref{grad_bound},  \ref{lem_g3} and Remark \ref{phi_def} that
 there exist rotations $R_h(x'):S\to SO(3)$ and $Q_h(x'):S\to SO(2)$ such that
  \begin{align*}
&\|\nabla_h y(x',x_3)-R_h(x')\|_{L^2(\Omega)}^2\le Ch^{\sigma-2},\;\|\nabla' R_h(x')\|_{L^2(S)}^2\le Ch^{\sigma-4},\\
&\|\int\limits_I\nabla' y'(x)dx_3-Q_h(x')\|_{L^2(S)}^2\le Ch^\sigma,\; \|\nabla' Q_h(x')\|_{L^2(S)}^2\le Ch^{\sigma-2}.
\end{align*}
We expand $W_i$ around the identity, $W_i(Id+A)=\frac{1}{2}\mathcal{Q}_3^i(A)+\eta_i(A)$, where $\mathcal{Q}_3^i(A)=\frac{\partial^2W_i(A)}{\partial A^2}(Id)(A,A)$ and $\eta_i(A)/|A|^2\to 0$ as $|A|\to 0$. If now $\omega_i(t)=\sup\limits_{|A|\le t}|\eta_i(A)|$ we have
 \begin{equation}
     \label{taylor}
     W_i(Id+A)\ge \frac12\Q_3^i(A)-\omega_i(|A|).
 \end{equation}
With the notation $r_h:=\frac{y_3^h}{(E_h)^{1/4}}$ it follows from \ref{A4}  that
  \begin{equation}
  \label{r_conv}
\nabla'r_h\otimes \nabla'r_h\rightharpoonup b(x) \quad\text{in} \;\; L^2(\Omega;\mathbb R^{2\times 2}).
   \end{equation}
   If $\sigma > 4$, then \eqref{y_i} yields 
   \begin{equation}
\nabla'r_h\to 0 \quad\text{in} \quad L^2(\Omega;\mathbb R^{2}),\label{r_conv1}
   \end{equation}
   and therefore, $\nabla'r_h\otimes \nabla'r_h \to 0$ in  $L^1(\Omega;\mathbb R^{2\times 2})$ and in view of \eqref{r_conv} 
    \begin{equation}
    \label{r_conv111}
\nabla'r_h\otimes \nabla'r_h\rightharpoonup 0 \quad\text{in} \quad L^2(\Omega;\mathbb R^{2\times 2}), \quad \text{ if } \sigma>4.
   \end{equation}
   We estimate
 \begin{eqnarray*}
      &&\liminf\limits_{h\to 0}\frac{1}{E_h}I^h(y^h)  
      \ge \liminf\limits_{h\to 0}
      (\frac{1}{E_h} \int\limits_\Omega W_1(\nabla'(y^h)'^TQ_hQ_h^T\nabla'(y^h)'+\nabla'y_3^h\otimes \nabla'y_3^h)dx +\frac{h^2}{E_h}\int\limits_\Omega W_2(L_h^T\nabla_hy^h)dx)  \\&=&\liminf\limits_{h\to 0}(\frac{h^2}{E_h}\int\limits_\Omega W_2(Id+\frac{\sqrt{E_h}}{h}G_h)dx \\ &&+
      \frac{1}{E_h} \int\limits_\Omega W_1((Id+\sqrt{E_h}F_h+\sqrt{E_h}K_h)^T(Id+\sqrt{E_h}F_h+\sqrt{E_h}K_h)+\sqrt{E_h}\nabla'r_h\otimes \nabla'r_h)dx) \\
      &=& \liminf\limits_{h\to 0}(\frac{h^2}{E_h}\int\limits_\Omega W_2(Id+\frac{\sqrt{E_h}}{h}G_h)dx +
      \frac{1}{E_h} \int\limits_\Omega W_1((Id+2\sqrt{E_h} \sym F_h +2\sym \sqrt{E_h}K_h   \\&&+E_hF_h^TF_h+ E_h K_h^TK_h+2E_h \sym (K_h^TF_h) +\sqrt{E_h}\nabla'r_h\otimes \nabla'r_h)dx).     
 \end{eqnarray*}
 Let us define $\chi_h$ as a characteristic function of the set
 \begin{eqnarray*}
 \Omega_h=\left\{\quad x\in \Omega: |G_h|\le 1/(E_h)^{1/8}, \; |K_h|\le 1/(E_h)^{1/8}, \quad |F_h|\le 1/(E_h)^{1/8},\; |\nabla 'r_h|\le 1/(E_h)^{1/8}\right\}.  
 \end{eqnarray*}
Then $\chi_h$ is bounded and $\chi_h\to 1$ in $L^1(\Omega)$, thus we have $\chi_hG_h\rightharpoonup G$, $\chi_hK_h\rightharpoonup K$, $\chi_hF_h\rightharpoonup F$, $\chi_hr_h\otimes r_h\rightharpoonup b(x')$ in $L^2(\Omega;\mathbb R^{2\times 2})$.
Therefore,
  \begin{eqnarray}
     \label{bel1}
   &&   \liminf\limits_{h\to 0}\frac{1}{E_h}I^h(y^h)\ge 
       \liminf\limits_{h\to 0}(\frac{h^2}{E_h}\int\limits_\Omega \chi_h W_2(Id+\frac{\sqrt{E_h}}{h}G_h)dx\nonumber\\&&+
      \frac{1}{E_h} \int\limits_\Omega \chi_h  W_1((Id+2\sqrt{E_h} \sym F_h+2\sym \sqrt{E_h}K_h\nonumber\\&&+E_hF_h^TF_h+E_h K_h^TK_h+2E_h \sym (K_h^TF_h)+\sqrt{E_h}\nabla'r_h\otimes \nabla'r_h)dx)\nonumber\\&&\ge \liminf\limits_{h\to 0}(\int\limits_\Omega (\chi_h\frac{1}{2} \Q_3^2(G_h)-\frac{h^2}{E_h} \chi_h\omega_2(\frac{\sqrt{E_h}}{h}G_h))dx\nonumber\\&&+
       \int\limits_\Omega (\chi_h  \frac{1}{2} \Q_3^1((2\sym F_h+2\sym K_h+\nabla'r_h\otimes \nabla'r_h)
       -C\chi_h E_h(|F_h|^4+|K_h|^4)-\frac{1}{E_h} \chi_h\omega_1(2\sqrt{E_h}\sym F_h\nonumber\\&&+2\sqrt{E_h}\sym K_h+E_hF_h^TF_h+E_h K_h^TK_h+2E_h\sym (K_h^TF_h)+\sqrt{E_h}\nabla'r_h\otimes \nabla'r_h))dx).
       \end{eqnarray}
       It is easy to see that whenever $\chi_h\ne 0$, $\frac{\sqrt{E_h}}{h}|G_h|\le \frac{E_h^{3/8}}{h}\to 0$ and
       \begin{equation}
\label{om1}
\frac{h^2}{E_h} \chi_h\omega_2\left(\left|\frac{\sqrt{E_h}}{h}G_h\right| \right)=\frac{\chi_h\omega_2(|\frac{\sqrt{E_h}}{h}G_h)|)}{|\frac{\sqrt{E_h}}{h}G_h|^2} |G_h|^2\to 0.
       \end{equation}
Moreover, 
  \begin{equation}
\label{om2}
E_h(|F_h|^4+|K_h|^4)\le \sqrt{E_h}\to 0.
        \end{equation}
If we denote 
\begin{eqnarray*}
U_h=2\sym F_h+2\sym K_h+\sqrt{E_h}F_h^TF_h+\sqrt{E_h}K_h^TK_h+2\sqrt{E_h}\sym K_h^TF_h+\nabla'r_h\otimes \nabla'r_h,
\end{eqnarray*}
         then
         \begin{eqnarray*}
|\sqrt{E_h}U_h|=|2\sqrt{E_h}\sym F_h+2\sqrt{E_h}\sym K_h+E_hF_h^TF_h+E_h K_h^TK_h\\+2E_h\sym K_h^TF_h+\sqrt{E_h}\nabla'r_h\otimes \nabla'r_h|\le C(E_h)^{1/4}\to 0
\end{eqnarray*}
     and    
         \begin{eqnarray*}
  \frac{1}{E_h} \chi_h\omega_1\left( \left|\sqrt{E_h}\sym F_h+2\sqrt{E_h}\sym K_h+E_hF_h^TF_h+E_h K_h^TK_h \right. \right. \\ \left. \left. +2E_h\sym K_h^TF_h+\sqrt{E_h}\nabla'r_h\otimes \nabla'r_h\right| \right) \\\le \chi_h\frac{\omega_1(|\sqrt{E_h}U_h|)}{|\sqrt{E_h}U_h|^2} (|2\sym F_h+2\sym K_h+\nabla'r_h\otimes \nabla'r_h|^2+ \sqrt{E_h}  )  \to 0,
         \end{eqnarray*}
which together with \eqref{bel1}--\eqref{om2} yields 
\eqref{inflim1} in the case $\sigma>4$. 
 \item[(ii)] 
We assume that $v$, $u$ and $\varphi$ are smooth and consider
\begin{eqnarray}
\label{recov}
\hat y^h(x',x_3):=\left(\begin{array}{c} x'\\hx_3   \end{array}\right)+\left(\begin{array}{c} h^{\sigma/2}u(x')\\h^{\sigma/2-1}v(x')   \end{array}  \right)+x_3 \left(\begin{array}{c} h^{\sigma/2}\varphi(x')\\ 1/2h^{\sigma/2}\EuScript L(\tilde G)\end{array}  \right),
\end{eqnarray}
 where $c=\EuScript L(\tilde G)$ is the element which realizes the minimum of $\Q_2^2$, i.e.
 $$\Q_2^2(\tilde G)=\Q_3^2(\tilde G+ce_3\otimes e_3)).$$
 \end{enumerate}
 Then 
 \begin{eqnarray*}
\nabla_h\hat y^h=\left(\begin{array}{cc} Id+h^{\sigma/2}\nabla'u+x_3h^{\sigma/2}\nabla'\varphi &   h^{\sigma/2-1}\varphi\\
 h^{\sigma/2-1}\nabla'v& 1+1/2h^{\sigma/2-1}\EuScript L(\tilde G)
\end{array}\right)+\left(\begin{array}{cc} o(h^{\sigma/2})   & o(h^{\sigma/2-1}) \\o(h^{\sigma/2-1})&o(h^{\sigma/2-1})   \end{array}\right)
\end{eqnarray*}
and
 \begin{eqnarray*}
(\nabla_h\hat y^h)^T \nabla_h\hat y^h=\left(\begin{array}{cc} Id+2h^{\sigma/2}\sym \nabla'u+2x_3h^{\sigma/2}\sym \nabla'\varphi &   h^{\sigma/2-1}(\varphi+\nabla'v)\\
 h^{\sigma/2-1}(\varphi+\nabla'v)& 1+h^{\sigma/2-1}\EuScript L(\tilde G))
\end{array}\right)\\+\left(\begin{array}{cc} o(h^{\sigma/2})   & o(h^{\sigma/2-1}) \\o(h^{\sigma/2-1})&o(h^{\sigma/2-1})   \end{array}\right).
\end{eqnarray*}
Finally, using the Taylor expansion we obtain
 \begin{eqnarray*}
\frac{1}{h^\sigma}I^h\to \frac{1}{2}\int\limits_\Omega  \Q_2^2(\tilde G)dx'+
       \frac{1}{2}\int\limits_S \Q_3^1(2\sym \nabla'u)dx'+\frac{1}{6}\int\limits_S \Q_3^1(\sym \nabla'\varphi(x'))dx'.
\end{eqnarray*}
For general $u,\varphi,v \in W^{1,2}(S)$ the assertion follows by suitable smooth approximations.
\end{proof}
Using Theorem \ref{th_conv1} and arguing as in \cite[Theorem 2]{FGM_Der} one can prove Theorem  \ref{th_j} .
\section{The case $\sigma=4$.}\label{sec:dimred2}
In case $\sigma=4$ there is a lack of compactness for the sequence $v^h$. Due to this fact, it seems impossible to perform a limit transition in the nonlinear term $\nabla'r_h\otimes \nabla'r_h$. This creates obstacles in the derivation of a $\Gamma$-limit via $v,u$ and $\varphi$. 
In this case we begin with the energy functionals containing second gradient terms  (see, e.g. \cite{Mueller})
\begin{eqnarray*}
I_2^h(y)=\int\limits_\Omega W_1(\nabla'y'(x)^T \nabla'y'(x)+\nabla' y_3\otimes  \nabla' y_3(x))dx+h^2\int\limits_{\Omega}W_2(\nabla_h y(x)) dx \\+lh^2 \int\limits_{\Omega} |\nabla_h^2 y(x)|^2 dx 
+c_1\int\limits_{\Omega} |\nabla' y_3(x)|^{4} dx.
\end{eqnarray*}
\begin{theorem}
      \label{th_conv2}
     Suppose that the assumptions of Theorem \ref{th_conv1} are satisfied. 
\begin{enumerate} 
\item[(i)](Compactness and lower bound) \\
 If 
 \begin{equation}
  \label{boun2}
\limsup\limits_{h\to 0}\frac{1}{h^4}I_2^h(y^h)<\infty,
  \end{equation}
  then properties \eqref{u_conve}--\eqref{phi_conve2}  are satisfied with $\sigma=4$.
 Moreover, 
\begin{align}
\label{str}
(v^h)(x')=\frac{1}{h} \int\limits_I y_3^hdx_3\rightharpoonup  v\quad \text{in} \;\; W^{2,2}(S),
\end{align}
 and
  \begin{eqnarray*}
     \liminf\limits_{h\to 0}\frac{1}{h^\sigma}I_2^h(y^h)\ge
 \left(\frac{1}{2}\int\limits_S \Q_2^2({\tilde G})dx'+\frac{1}{6}\int\limits_S \Q_3^1(\sym \nabla'\varphi(x'))dx_1\right.\\\left.+
    \frac{1}{2}\int\limits_S \Q_3^1((2\sym \nabla' u(x')+\nabla'v\otimes \nabla'v(x'))+l\int\limits_S|\nabla'^2v(x')|^2dx'\right.\\\left.+l\int\limits_S|\nabla '\varphi(x')|^2dx'+c_1\int\limits_{S} |\nabla' v(x')|^{4} dx'\right).
 \end{eqnarray*}    
 \item[(ii)](Optimality of lower bound.)\\
 If  $u,\varphi \in W^{1,2}(S)$, $v \in W^{2,2}(S)$ then there exists $\hat y^h$ such that \eqref{y_i} and \eqref{u_conve}--\eqref{phi_conve2} hold and 
 \begin{eqnarray*}
      \liminf\limits_{h\to 0}\frac{1}{E_h}I^h(\hat y^h)= 
\left(\frac{1}{2}\int\limits_S \Q_2^2({\tilde G})dx'+\frac{1}{6}\int\limits_S \Q_3^1(\sym \nabla'\varphi(x'))dx_1\right.\\\left.+
    \frac{1}{2}\int\limits_S \Q_3^1(2\sym \nabla' u(x')+\nabla'v\otimes \nabla'v(x'))+l\int\limits_S|\nabla'^2v(x')|^2dx'\right.\\\left.+l\int\limits_S|\nabla'\varphi(x')|^2dx'+c_1\int\limits_{S} |\nabla' v(x')|^{4} dx'\right). 
 \end{eqnarray*}    
 \end{enumerate}
  \end{theorem}
  \begin{proof}
Arguing as in Theorem \ref{th_conv1} one can  show \eqref{y_i} and \eqref{u_conve}--\eqref{phi_conve2}.
It follows from \eqref{boun2} that
\begin{equation*}
 \int\limits_{\Omega}|\nabla_h^2 y(x)|^2 dx \le Ch^2.
\end{equation*}
In particular this means for $i=1,2$
\begin{align}
					&\int\limits_{\Omega}|\frac{1}{h}\nabla'^2 y'(x)|^2 dx \le C,\;&\int\limits_{\Omega}|\frac{1}{h}\nabla'^2 y_3(x)|^2 dx \le C,\label{l1}\\
					&\int\limits_{\Omega}|\frac{1}{h^2}\partial_i\partial_3 y'(x)|^2 dx \le C,\;\;&\int\limits_{\Omega}|\frac{1}{h^2}\partial_i\partial_3 y_3(x)|^2 dx \le C,\label{l2}\\
					&\int\limits_{\Omega}|\frac{1}{h^3}\partial_3^2 y'(x)|^2 dx \le C,\;
					&\int\limits_{\Omega}|\frac{1}{h^3}\partial_3^2 y_3(x)|^2 dx \le C.	\label{l3}	
				\end{align}
It follows from the second estimate in \eqref{l2} and the first estimate in \eqref{l3} that 
\begin{align}
  \frac{1}{h}\partial_i\partial_3 y_3(x) \to 0\;\text{in}\;L^2(\Omega)\text{\qquad and}\label{l10}\\
  \frac{1}{h^2}\partial_3^2 y'(x)\to 0\;\;\text{in}\;L^2(\Omega).\label{l11}
\end{align}
Taking into account \eqref{l11} and arguing as in Theorem \ref{th_conv1} we get \eqref{phi_conve1} which together with \eqref{p_conv2} yields
\begin{eqnarray*}
\frac{1}{h^2}\partial_3 y'(x)=\frac{1}{h^2}(I-Q_h^T)\partial_3 y'(x)+\frac{1}{h^2}Q_h^T\partial_3 y'(x)\rightharpoonup \varphi(x') \;\;\;\text{in}\;\;L^1(\Omega).
\end{eqnarray*}
Consequently, 
\begin{equation*}
\frac{1}{h^2}\nabla'\partial_3 y'(x)\rightharpoonup \nabla'\varphi(x')\;\;\;\text{in}\;W^{-1,1}(\Omega).
\end{equation*}
Therefore, \eqref{w1} together with the first estimate in \eqref{l2} we obtain
\begin{equation}
\label{w1}
\frac{1}{h^2}\nabla'\partial_3 y'(x)\rightharpoonup \nabla'\varphi(x')\;\;\;\text{in}\;L^{2}(\Omega).
\end{equation}
Moreover, from \eqref{l10} we have
\begin{equation}
\label{c2}
\|\frac{1}{h}\nabla'y_3^h-\frac{1}{h}\int\limits_I\nabla'y_3^hdx_3\|_{L^2(\Omega)}\le C\|\frac{1}{h}\nabla'\partial_3y_3^h\|_{L^2(\Omega)}\to 0.
\end{equation}
Consequently, \eqref{v_conve} and \eqref{c1} yield
\begin{equation*}
\frac{1}{h}\nabla'^2y_3^h\rightharpoonup \nabla'^2v(x')\;\;\;\text{in}\;W^{-1,2}(\Omega),
\end{equation*}
which together with the second estimate in \eqref{l1} leads to
\begin{equation}
\label{c1}
\frac{1}{h}\nabla'^2y_3^h\rightharpoonup \nabla'^2v(x')\;\;\;\text{in}\;L^{2}(\Omega).
\end{equation}
Next, \eqref{u_conve}, \eqref{phi_conve2},  \eqref{p_conv2} imply 
\begin{equation*}
\frac{1}{h}(\nabla'y'^h-Id)\rightharpoonup 0\;\;\;\text{in}\;L^{1}(\Omega).
\end{equation*}
Consequently, it follows from the first estimate in \eqref{l1} that
\begin{equation}
\label{c3}
\frac{1}{h}\nabla'^2y'^h\rightharpoonup 0\;\;\;\text{in}\;L^{2}(\Omega).
\end{equation}
Convergence \eqref{str} follows from \eqref{c1}. This yields $\frac{1}{h^2}\nabla'y_3^h\otimes \nabla'y_3^h\to \nabla'v'\otimes \nabla'v'$ in $L^2(\Omega)$.  Then we  choose as a recovery sequence \eqref{recov} and take into consideration that the second term in \eqref{l2} and terms in \eqref{l3} are equal to zero on the recovery sequence. Arguing as in Theorem \ref{th_conv1} and using \eqref{w1}, \eqref{c1}, \eqref{c3} we conclude the statement of the theorem. 
  \end{proof}
In this case the limiting functional contains a superfluous last term  which is not a part of classical nonlinear Reissner-Mindlin model. It  is an open question whether the nonlinear model can be derived without this term. In general, the nonlinear Reissner-Mindlin model lacks compactness, i.e. the nonlinearity is supercritical. 
\subsubsection*{Acknowledgements}
The first author was supported by Einstein foundation (EGR-2022-731). This work benefitted from discussions with K. Buryachenko, A. Glitzky and M. Liero in the context of project AA2-21 within MATH+ which is funded by the Deutsche Forschungsgemeinschaft (DFG, German Research Foundation) under Germany´s Excellence Strategy – The Berlin Mathematics Research Center MATH+ (EXC-2046/1, project ID 390685689).


\vskip-10mm
\begin{thebibliography}{99}

\bibitem{ABP} Acerbi, E., Buttazzo, G. and Percivale, D., 
 A variational definition of the strain energy for an elastic string, J. Elasticity, 25(2), 137--148 (1991). https://doi.org/10.1007/BF00042462
 
\bibitem{ADS}Agostiniani, V. and DeSimone, A., Rigorous derivation of active plate models for thin sheets of nematic elastomers, Math. Mech. Solids, 25(10), 1804--1830 (2020). https://doi.org/10.1177/1081286517699991

 \bibitem{BSN} Balzani, D., Schröder, J. and Neff, P., Applications of anisotropic polyconvex energies: thin shells and biomechanics of arterial walls. In: Schröder, J., Neff, P. (eds) Poly-, Quasi- and Rank-One Convexity in Applied Mechanics. CISM International Centre for Mechanical Sciences, vol 516. Springer, Vienna, (2010). https://doi.org/10.1007/978-3-7091-0174-2\_5

\bibitem{BGNPP} Bartels, S., Griehl, M., Neukamm, S., Padilla-Garza, D. and Palus, C., A nonlinear bending theory for nematic LCE plates, Math. Models Methods Appl. Sci., 33(07), 1437-1516 (2023). https://doi.org/10.1142/S0218202523500331

\bibitem{DBS}de Benito Delgado, M. and Schmidt, B., A hierarchy of multilayered plate models, ESAIM Control Optim. Calc. Var. , 27, S16 (2021). https://doi.org/10.1051/cocv/2020067

\bibitem{BB} Bonet J. and  Burton A.J., A simple orthotropic, transversely isotropic hyperelastic constitutive
equation for large strain computations, Comput. Methods Appl. Mech. Engrg. 162, 151--164 (1998). https://doi.org/10.1016/S0045-7825(97)00339-3

\bibitem{braides}Braides, A., Gamma-convergence for Beginners (Vol. 22), Clarendon Press (2002). https://doi.org/10.1093/acprof:oso/9780198507840.001.0001

\bibitem{conti-zwicknagl}
Conti, S. and Zwicknagl, B., The tapering length of needles in martensite/martensite macrotwins, Arch. Rational Mech. Anal., 247 (4), 63 (2023). https://doi.org/10.1007/s00205-023-01882-9



\bibitem{DalMaso} Dal Maso, G., An introduction to $\Gamma$-convergence. Birhäuser, 1993.

\bibitem{Davoli} Davoli, E., Thin-walled beams with a cross-section of arbitrary geometry: derivation of linear theories starting from 3D nonlinear elasticity, Adv. Calc. Var., 6(1), 33--91 (2013).  https://doi.org/10.1515/acv-2011-0003 



\bibitem{EK}Engl D. and Kreisbeck, C., Theories for incompressible rods: A rigorous derivation via $\Gamma$-convergence, Asymptot. Anal., 124(1-2), 1--28 (2020).  https://doi.org/10.3233/ASY-201636

  \bibitem{FPP} Falach L., Paroni, R. and  Podio-Guidugli, P., A justification of the Timoshenko beam model through $\Gamma$-convergence, Anal. Appl., 15, No. 2,  261-277 (2017). https://doi.org/10.1142/S0219530515500207

\bibitem{FJM}
  Friesecke, G., James and R.D., Müller, S., A theorem on geometric rigidity and the derivation of nonlinear plate theory from three-dimensional elasticity, Comm. Pure Appl. Math., 55(11), 1461--1506 (2002).  https://doi.org/10.1002/cpa.10048

 \bibitem{FGM_Der}  Friesecke, G., James, R.D., Müller, S., A hierarchy of plate models derived from nonlinear elasticity by Gamma-convergence, Arch. Rational Mech. Anal., 180, 183--236 (2006). https://doi.org/10.1007/s00205-005-0400-7

 \bibitem{FGM_Rig}  Friesecke, G., James, R.D., Müller, S., Rigorous derivation of nonlinear plate theory and geometric rigidity, C.R. Acad. Sci.Paris. S\'er I, 334, 173--178 (2002).  https://doi.org/10.1016/S1631-073X(02)02133-7
  
\bibitem{FGMM} Friesecke, G., James, R.D., Müller, S., Mora, M.G., Derivation of nonlinear bending theory for shells from three dimensional nonlinear elasticity by Gamma-convergence, C.R. Acad. Sci.Paris. S\'er I, 336, 697--702 (2003).  https://doi.org/10.1016/S1631-073X(03)00028-1

\bibitem{GOW}Griso, G., Orlik, J. and Wackerle, S., Asymptotic behavior for textiles in von-Kármán regime, J. Math. Pures Appl., 144, 164-193  (2020). https://doi.org/10.1016/j.matpur.2020.10.002

\bibitem{GG}
 Ginster J. and  Gladbach P.,
The Euler-Bernoulli limit of thin brittle linearized elastic beams, 
J. Elast., 156, (2023). https://doi.org/10.1007/s10659-023-10040-x

\bibitem{LDR1} Le Dret, H. and Raoult, A., The nonlinear membrane model as variational limit of nonlinear three-dimensional elasticity, J. Math. Pures Appl., 74(6), 549--578 (1995).

\bibitem{LDR2} Le Dret, H. and Raoult, A., The membrane shell model in nonlinear elasticity: a variational asymptotic derivation, J. Nonlinear Sci., 6(1), 59--84 (1996). https://doi.org/10.1007/BF02433810

\bibitem{LM}  Lecumberry, M. and  Müller, S., Stability of slender bodies under compression and validity of the von Karman theory, Arch. Rational Mech. Anal., 193, 255--310 (2009). https://doi.org/10.1007/s00205-009-0232-y

\bibitem{lewicka} Lewicka, M., Calculus of Variations on Thin Prestressed Films, Springer International Publishing (2023).

  \bibitem{M}   Mindlin, R. D.,  Influence of rotatory inertia and shear on flexural motions of isotropic, elastic plates, J. Appl. Mech., 18, 31--38 (1951). https://doi.org/10.1115/1.4010217


 \bibitem{Mueller} Müller, S.,  Variational models for microstructure and phase transitions. In: Hildebrandt, S., Struwe, M. (eds) Calculus of Variations and Geometric Evolution Problems. Lecture Notes in Mathematics, vol 1713. Springer, Berlin, Heidelberg (1999). https://doi.org/10.1007/BFb0092670

\bibitem{MM1} Mora, M.G., and  Müller, S., Derivation of the nonlinear bending-torsion theory for inextensible rods by $\Gamma $-convergence, Calc. Var. Partial Differential Equations, 18(3), (2003). https://doi.org/10.1007/s00526-003-0204-2

\bibitem{MM2} Mora, M.G., and Müller, S., A nonlinear model for inextensible rods as a low energy $\Gamma$-limit of three-dimensional nonlinear elasticity, Ann. Inst. H. Poincar\'e C Anal. Non Lin\'eaire, 21(3), 271-293 (2004). https://doi.org/10.1016/j.anihpc.2003.08.001

\bibitem{mueller} Müller, S., Mathematical problems in thin elastic sheets: scaling limits, packing, crumpling and singularities, In Vector-Valued Partial Differential Equations and Applications: Cetraro, Italy 2013 (pp. 125-193). Cham: Springer International Publishing (2017). https://doi.org/10.1007/978-3-319-54514-1\_3

\bibitem{nhj}Neff, P., Hong, K. I. and Jeong, J., The Reissner–Mindlin plate is the $\Gamma$-limit of Cosserat elasticity, Math. Models Methods Appl. Sci., 20(09), 1553-1590 (2010). https://doi.org/10.1142/S0218202510004763

\bibitem{Neu}Neukamm, S., Rigorous Derivation of a Homogenized Bending-Torsion Theory for Inextensible Rods from Three-Dimensional Elasticity, Arch. Rational Mech. Anal., 206(2), 645--706 (2012). https://doi.org/10.1007/s00205-012-0539-y


\bibitem{NV} Neukamm, S. and Vel\v{c}i\'{c}, I., Derivation of a homogenized von-K\'{a}rm\'{a}n plate theory from 3D nonlinear elasticity, Math. Models Methods Appl. Sci., 23(14), 2701-2748 (2013). https://doi.org/10.1142/S0218202513500449

 \bibitem{Ogden} Ogden, R.W., Nonlinear Continuum Mechanics and Modeling the Elasticity of Soft Biological Tissues with a Focus on Artery Walls. In: Holzapfel, G., Ogden, R. (eds) Biomechanics: Trends in Modeling and Simulation. Studies in Mechanobiology, Tissue Engineering and Biomaterials, vol 20. Springer, Cham. (2017). https://doi.org/10.1007/978-3-319-41475-1\_3

 \bibitem{Pan} Pantz, O., On the justification of the nonlinear inextensional plate model, Arch. Rational Mech. Anal. 167(3), 179-209 (2003). https://doi.org/10.1007/s00205-002-0238-1

\bibitem{PPG} Paroni, R. and  Podio-Guidugli, P., On variational dimension reduction in structure mechanics, J. Elasticity, 118(1), 1--13 (2015). https://doi.org/10.1007/s10659-014-9473-6
 
 \bibitem{PPT}  Paroni, R.,  Podio-Guidugli, P. and Tomasetti G. The Reissner-Mindlin plate theory via $\Gamma$-convergence, C.R. Acad. Sci. Paris, Ser. I, 343, 437--440 (2006).   https://doi.org/10.1016/j.crma.2006.08.006

 \bibitem{PPT2} Paroni, R., Podio-Guidugli, P. and Tomassetti, G., A justification of the Reissner–Mindlin plate theory through variational convergence, Anal. Appl., 5(02), 165-182 (2007). https://doi.org/10.1142/S0219530507000936

\bibitem{R}    Reissner, E.,  The effect of transverse shear deformation on the bending of elastic plates, J. Appl. Mech., Vol. 12, 68–77  (1945). https://doi.org/10.1115/1.4009435

\bibitem{Scardia}Scardia, L., Asymptotic models for curved rods derived from nonlinear elasticity by $\Gamma$-convergence, Proc. Roy. Soc. Edinburgh Sect. A, 139(5), 1037-1070 (2009). https://doi.org/10.1017/S0308210507000194

\bibitem{Scardia2}Scardia, L., The nonlinear bending–torsion theory for curved rods as $\Gamma$-limit of three-dimensional elasticity, Asymptot. Anal., 47(3-4), 317--343 (2006). https://doi.org/10.3233/ASY-2006-760

\bibitem{BS} Schmidt, B., A Griffith--Euler--Bernoulli theory for thin brittle beams derived from nonlinear models in variational fracture mechanics,
  Math. Models Methods Appl. Sci.,
  27 (09), 1685--1726 (2017). https://doi.org/10.1142/S0218202517500294


  \bibitem{SN} Schröder, J., Neff, P., Invariant formulation of hyperelastic transverse isotropy
based on polyconvex free energy functions,  Int. J. Solids Struct., Vol. 40  401--445 (2003). https://doi.org/10.1016/S0020-7683(02)00458-4

	 	\end{thebibliography}
\end{document}